%
\def\conv{\mathop{\vrule height2,6pt depth-2,3pt 
    width 5pt\kern-1pt\rightharpoonup}}

\advance\vsize by 1 true cm
%
\def\dess #1 by #2 (#3){
  \vbox to #2{
    \hrule width #1 height 0pt depth 0pt
    \vfill
    \special{picture #3} 
    }
  }

\def\dessin #1 by #2 (#3 scaled #4){{
  \dimen0=#1 \dimen1=#2
  \divide\dimen0 by 1000 \multiply\dimen0 by #4
  \divide\dimen1 by 1000 \multiply\dimen1 by #4
  \dess \dimen0 by \dimen1 (#3 scaled #4)}
  }
%
\def \trait (#1) (#2) (#3){\vrule width #1pt height #2pt depth #3pt}
\def \fin{\hfill
	\trait (0.1) (5) (0)
	\trait (5) (0.1) (0)
	\kern-5pt
	\trait (5) (5) (-4.9)
	\trait (0.1) (5) (0)
\medskip}
%


\font\sevenbf=cmbx7

\baselineskip=15pt
\abovedisplayskip=15pt plus 4pt minus 9pt
\belowdisplayskip=15pt plus 4pt minus 9pt
\abovedisplayshortskip=3pt plus 4pt
\belowdisplayshortskip=9pt plus 4pt minus 4pt
\let\epsilon=\varepsilon

\def\biblio #1 #2\par{\parindent=30pt\item{}\kern -30pt\rlap{[#1]}\kern
30pt #2\smallskip}
 %
\catcode`\@=11
\def\@lign{\tabskip=0pt\everycr={}}
\def\equations#1{\vcenter{\openup1\jot\displ@y\halign{\hfill\hbox
{$\@lign\displaystyle##$}\hfill\crcr
#1\crcr}}}
\catcode`\@=12
%
\def\pmb#1{\setbox0=\hbox{#1}%
\hbox{\kern-.04em\copy0\kern-\wd0
\kern.08em\copy0\kern-\wd0
\kern-.02em\copy0\kern-\wd0
\kern-.02em\copy0\kern-\wd0
\kern-.02em\box0\kern-\wd0
\kern.02em}}
%
\def\undertilde#1{\setbox0=\hbox{$#1$}
\setbox1=\hbox to \wd0{$\hss\mathchar"0365\hss$}\ht1=0pt\dp1=0pt
\lower\dp0\vbox{\copy0\nointerlineskip\hbox{\lower8pt\copy1}}}
%

%

\def\maj#1#2,{\rm #1\sevenrm #2\rm{}}
\def\Maj#1#2,{\bf #1\sevenbf #2\rm{}}
\outer\def\lemme#1#2 #3. #4\par{\medbreak
\noindent\maj{#1}{#2},\ #3.\enspace{\sl#4}\par
\ifdim\lastskip<\medskipamount\removelastskip\penalty55\medskip\fi}

\def\Remark #1. {\noindent{\Maj REMARK,\ \bf #1. }}

\outer\def\Lemme#1#2 #3. #4\par{\medbreak
\noindent\Maj{#1}{#2},\ \bf #3.\rm\enspace{\sl#4}\par
\ifdim\lastskip<\medskipamount\removelastskip\penalty55\medskip\fi}



\def\Notation #1. {\noindent{\Maj NOTATION,\ \bf #1. }}

\def\Example #1. {\noindent{\Maj EXAMPLE,\ \bf #1. }}

\hfuzz=1cm


\catcode`\ˆ=\active     \def ˆ{\`a}
\catcode`\‰=\active     \def ‰{\^a}
\catcode`\=\active     \def {\c c}
\catcode`\Ž=\active    \def Ž{\'e} 

\catcode`\=\active   \def {\`e}
\catcode`\=\active   \def {\^e}
\catcode`\'=\active   \def '{\"e}
\catcode`\"=\active   \def "{\^\i}
\catcode`\•=\active   \def •{\"\i}
\catcode`\™=\active   \def ™{\^o}
\catcode`\š=\active   \defš{}
\catcode`\=\active   \def {\`u}
\catcode`\ž=\active   \def ž{\^u}
\catcode`\Ÿ=\active   \def Ÿ{\"u}
\catcode`\ =\active   \def  {\tau}
\catcode`\¡=\active   \def ¡{\circ}
\catcode`\¢=\active   \def ¢{\Gamma}
\catcode`\¤=\active   \def ¤{\S\kern 2pt}
\catcode`\¥=\active   \def ¥{\puce}
\catcode`\§=\active   \def §{\beta}
\catcode`\¨=\active   \def ¨{\rho}
\catcode`\©=\active   \def ©{\gamma}
\catcode`\­=\active   \def ­{\neq}
\catcode`\°=\active   \def °{\ifmmode\ldots\else\dots\fi}
\catcode`\±=\active   \def ±{\pm}
\catcode`\²=\active   \def ²{\le}
\catcode`\³=\active   \def ³{\ge}
\catcode`\µ=\active   \def µ{\mu}
\catcode`\¶=\active   \def ¶{\delta}
\catcode`\·=\active   \def ·{\Sigma}
\catcode`\¸=\active   \def ¸{\Pi}
\catcode`\¹=\active   \def ¹{\pi}
\catcode`\»=\active   \def »{\Upsilon}
\catcode`\¾=\active   \def ¾{\alpha}
\catcode`\À=\active   \def À{\cdots}
\catcode`\Â=\active   \def Â{\lambda}
\catcode`\Ã=\active   \def Ã{\sqrt}
\catcode`\Ä=\active   \def Ä{\varphi}
\catcode`\Å=\active   \def Å{\xi}
\catcode`\Æ=\active   \def Æ{\Delta}
\catcode`\Ç=\active   \def Ç{\cup}
\catcode`\È=\active   \def È{\cap}
\catcode`\Ï=\active   \def Ï{\oe}
\catcode`\Ñ=\active   \def Ñ{\to}
\catcode`\Ò=\active   \def Ò{\in}
\catcode`\Ô=\active   \def Ô{\subset}
\catcode`\Õ=\active   \def Õ{\superset}
\catcode`\Ö=\active   \def Ö{\over}
\catcode`\×=\active   \def ×{\nu}
\catcode`\Ù=\active   \def Ù{\Psi}
\catcode`\Ú=\active   \def Ú{\Xi}
\catcode`\Ü=\active   \def Ü{\omega}
\catcode`\Ý=\active   \def Ý{\Omega}
\catcode`\ß=\active   \def ß{\equiv}
\catcode`\à=\active   \def à{\chi}
\catcode`\á=\active   \def á{\Phi}
\catcode`\ä=\active   \def ä{\infty}
\catcode`\å=\active   \def å{\zeta}
\catcode`\æ=\active   \def æ{\varepsilon}
\catcode`\è=\active   \def è{\Lambda}  
\catcode`\é=\active   \def é{\kappa}
\catcode`\ë=\active   \defë{\Theta}
\catcode`\ì=\active   \defì{\eta}
\catcode`\í=\active   \defí{\theta}
\catcode`\î=\active   \defî{\times}
\catcode`\ñ=\active   \defñ{\sigma}
\catcode`\ò=\active   \defò{\psi}

\def\date{\number\day\
\ifcase\month \or janvier \or f\'evrier \or mars \or avril \or mai \or juin \or juillet \or ao\^ut  \or
septembre \or octobre \or novembre \or d\'ecembre \fi
\ \number\year}

\font \ggras=cmb10 at 11pt

\def\sym{\fam\comfam\com}
\font\tensym=msbm10
\font\sevensym=msbm7
\font\fivesym=msbm5
\newfam\symfam
\textfont\symfam=\tensym
\scriptfont\symfam=\sevensym
\scriptfont\symfam=\fivesym
\def\sym{\fam\symfam\relax}

\def\R{{\sym R}}

\def\Z{{\sym Z}}


\font \Gggras=cmb10 at 18pt

\font \ggras=cmb10 at 16pt

\def\biblio #1 #2\par{\parindent=30pt\item{}\kern -30pt\rlap{[#1]}\kern 30pt #2\smallskip}
\def\biblio #1 #2\par{\parindent=30pt\item{[]}\kern -30pt\rlap{[#1]}\kern 30pt #2\smallskip}

\vskip 0.5cm
\centerline{\Gggras Interior error estimate for periodic  homogenization}
\vskip 5mm
\centerline{G. Griso}
\vskip 8mm
\centerline{  Laboratoire J.-L. Lions--CNRS, Bo\^\i te courrier 187, Universit\'e  Pierre et
Marie Curie, }
\centerline{ 4~place Jussieu, 75005 Paris, France,  Email: griso@ann.jussieu.fr}
\vskip 8mm
\noindent{\ggras \bf Abstract. } {\sevenrm 

In a previous article about the homogenization of the classical problem of diffusion in a bounded domain  with sufficiently
smooth boundary   we proved that the error is of order $\scriptstyle\varepsilon^{1/2}$. Now, for an open set $\scriptstyle\Omega$ with
sufficiently smooth boundary ($\scriptstyle {\cal C}^{1,1}$) and  homogeneous Dirichlet or Neuman limits conditions we show 
that in any open set strongly included in $\scriptstyle \Omega$ the error is of order $\scriptstyle
\varepsilon$. If the open set $\scriptstyle \Omega\subset\R^n$ is of polygonal ($\scriptstyle n=2$) or polyhedral $\scriptstyle
(n=3)$ boundary  we also give the global and interrior error estimates. }
\bigskip
\noindent{\ggras \bf R\'esum\'e. } {\sevenrm  Nous avons dŽmontrŽ dans un prŽcŽdent article sur l'homogŽnŽisation du
problme type de la diffusion dans un domaine bornŽ de frontire rŽgulire que l'erreur est d'ordre $\scriptstyle
\varepsilon^{1/2}$. On montre maintenant pour un ouvert $\scriptstyle \Omega$  de fronti\`ere  r\'eguli\`ere
($\scriptstyle {\cal C}^{1,1}$) avec les conditions aux limites homog\`enes de Dirichlet ou de Neumann  que dans tout 
ouvert fortement inclus dans $\scriptstyle\Omega$ l'erreur est de l'ordre de $\scriptstyle \varepsilon$.
Si l'ouvert $\scriptstyle \Omega\subset \R^n$ est de fronti\`ere polygonale ($\scriptstyle n=2$) ou poly\'edrale $\scriptstyle
(n=3)$ on donne Žgalement les estimations globale  et intŽrieure de l'erreur. }
\bigskip
\noindent{\bf 1. Introduction}
\vskip 1mm
This paper follows two previous studies [4,5] of the error estimates in the classical periodic homogenization problem. 
 The first error estimates in periodic homogenization problem have been given by  Bensoussan, Lions and Papanicolaou [1], by
Oleinik, Shamaev and Yosifian  [7], and by Cioranescu and Donato [3].  In all these works, the result  is proved under the
assumption that the correctors belong to $W^{1,\infty}(Y)$, $Y=]0,1[^n$ being the reference cell. The  estimate is of order  
$\varepsilon^{1/ 2}$. The  additional  regularity of the correctors holds true  when the coefficients of the operator are very
regular  which is not necessarily the situation in homogenization. In [4]  we  obtained an error estimate without any regularity
hypothesis on the correctors but we supposed that the solution of the homogenized problem belonged to
$W^{2,p}(\Omega)$ ($p>n$). The exponent of $\varepsilon$ in the error estimate is inferior to $1/2$ and depends on $n$ and
$p$.   In [5]  we  obtained an error estimate without any regularity hypothesis on the correctors but we supposed that the solution
of the homogenized problem belonged to $H^2(\Omega)$. This holds true with a smooth boundary   and  homogeneous Dirichlet or
Neuman limits conditions.  The exponent of $\varepsilon$ in the error estimate is equal to $1/2$. 
\smallskip
The aim of this work is to give the interior error estimate and new error estimate  with minimal hypothesis on the boundary of
$\Omega$. 
\smallskip
The paper is organized as follows. Section 2 is dedicated to some  projection theorems. Among them Theorems 2.3 and 2.6
are essential tools  to obtain new estimates. These theorems are related to  the periodic unfolding method (see [2] and [5]).  
We show that for any $\phi$ in $H^1(\Omega)$, where $\Omega$ is a bounded open  set of $\R^n$  with Lipschitz boundary,
there exists a function $\widehat{\phi}_\varepsilon$ in $ L^2(\Omega; H^1_{per}(Y ))$, such that the distance between the
unfolded
${\cal T}_\varepsilon(\nabla_x\phi)$ and $\nabla_x\phi+\nabla_y\widehat{\phi}_\varepsilon$ is of order 
$\varepsilon$ in the space $[L^2(Y ; (H^1(\Omega))^{'})]^n$ (Theorem 2.3) and  is of order  
$\varepsilon^s$ in the space $[L^2(Y ; (H^s(\Omega))^{'})]^n$, $0<s<1$, (Theorem 2.6), provided that the norm of
gradient $\phi$  in a neighbourhood (of thickness $4\varepsilon\sqrt n$) of the boundary of $\Omega$ is less than 
$\varepsilon^{1/2}$  in the first case and less than  $\varepsilon^{s/2}$  in the second case. 
\smallskip
 In  Theorem 3.2 in Section 3.1,  we suppose  that  $\Omega$ has a smooth boundary, that the right handside of the
homogenization problem belongs to  $L^2(\Omega)$  and  we consider the  homogeneous Dirichlet or Neumann limits conditions.
By transposition and thanks to Theorem 2.3 we show  that the $L^2$ error estimate is of order $\varepsilon$ and then we obtain
the interior error estimate of the same order.  The required condition in Theorem 2.3 is obtained thanks to the estimates of
Theorems 4.1 and 4.2 of [5]. 
\smallskip
In Theorem 3.3  in Section 3.2,   we suppose that the domain   $\Omega$ is of polygonal ($n=2$) or polyhedral 
$(n=3)$ boundary and the right handside of the homogenization problem in $L^2(\Omega)$. We show that the $H^1$ error
estimate is at the most of order $\varepsilon^{1/4}$ and that the $L^2$ and the interior  error estimates  are at the most of order
$\varepsilon^{1/2}$.
\smallskip
We  use   the notation of [2] and [5]  throughout this  study. In this article, the constants appearing in the estimates are
independent from $\varepsilon$.
\bigskip
\noindent{\bf 2. Preliminary results}
\medskip
\noindent Let  $\Omega$ be a  bounded domain in $\R^n$ with lipchitzian boundary. We put
$$\eqalign{ 
&\widehat{\Omega}_{\varepsilon, k}=\Bigl\{x\in \R^n\; |\; dist(x,\partial\Omega)<k\sqrt n \varepsilon\bigr\},\qquad
\widetilde{\Omega}_{\varepsilon, k}=\Bigl\{x\in \R^n\; |\; dist(x, \Omega)<k\sqrt n \varepsilon\bigr\},
\qquad k\in\{1,\,2,\, 3,\, 4\},\cr
&\Omega_\varepsilon =\hbox{interior}\Bigl(\bigcup_{\xi\in\Xi_\varepsilon}\varepsilon(\xi+\overline{Y}\Bigr),\qquad
\Xi_\varepsilon=\bigl\{\xi\in \Z^n\; |\; \varepsilon(\xi+ Y)\cap \Omega\not=\emptyset\bigr\},\qquad Y=]0,1[^n,\cr}$$ 
\noindent where the open set $Y=]0,1[^n$ is the reference cell and where  $\varepsilon$ is a strictly positive real.  We have
$$\Omega\i \Omega_\varepsilon\in\widetilde{\Omega}_{\varepsilon, 1}$$ For almost any $x\in \R^n$, there exists a unique
element in $\Z^n$ denoted  $[x]$  such that 
$$x=[x]+\{x\},\qquad \{x\}\in Y.$$ The running point of $\Omega$ is denoted  $x$, and the
running point of $Y$ is denoted  $y$.
\medskip
\noindent{\bf 2.1 Projection theorems in  $L^2(Y; (H^1(\Omega))^{'})$.}
\medskip
\noindent{\bf Lemma 2.1 : }{\it There exists a linear and continuous extension operator ${\cal P}_\varepsilon$ from
$H^1(\Omega)$  into $H^1(\widetilde{\Omega}_{\varepsilon, 3})$ such that
$$\left\{\eqalign{ 
&||\nabla{\cal P}_\varepsilon(\phi)||_{[L^2(\widetilde{\Omega}_{\varepsilon, 3})]^n}\le C
||\nabla\phi||_{[L^2(\Omega)]^n}\qquad ||\nabla{\cal P}_\varepsilon(\phi)||_{[L^2(\widetilde{\Omega}_{\varepsilon, 3}
\setminus\Omega)]^n}\le C ||\nabla\phi||_{[L^2(\Omega\setminus\widehat{\Omega}_{\varepsilon, 3})]^n}\cr
&|| {\cal P}_\varepsilon(\phi)||_{L^2(\widetilde{\Omega}_{\varepsilon, 3})}\le
C\bigl\{||\phi||_{L^2(\Omega)}+\varepsilon ||\nabla\phi||_{[L^2(\Omega\setminus \widehat{\Omega}_{\varepsilon,
3})]^n}\bigr\}\cr}\right.\leqno(2.1)$$ \noindent  The constants  depend only on $n$ and  $\partial\Omega$.}

\noindent{\bf Proof : } There exists a finite open covering $(\Omega_j)_j$ of the boundary $\partial
\Omega$ such that for each $j$ there exists a Lipschitz diffeomorphism $\theta_j$ which maps $\Omega_j$ to the open set
${\cal O}=]-1, 1[^{n-1}\times ]-1, 1[$ of $\R^n$ and $\Omega_j\cap \Omega$ to the open set ${\cal O}_+=]-1, 1[^{n-1}\times
]0, 1[$. To the covering of $\partial\Omega$ we associate a partition of the unity 
$$\phi_j\in {\cal C}^1_0(\Omega_j),\qquad \sum_j\phi_j=1\qquad \hbox{ in a neighbourhood of}\quad
\partial\Omega.$$ Let
 $\psi$ be in $H^1(\Omega)$. The function $(\phi_j\psi)\circ\theta_j^{-1}$ belongs to $H^1({\cal O}_+)$. We  use a reflexion
argument to extend this function to an element $\widetilde{\psi}_j$ belonging to $H^1({\cal O})$. In the neighbourhood of the
boundary of $\Omega$ the extension is equal to $\displaystyle\sum_j\widetilde{\psi}_j\circ
\theta_j$. This immediately  gives the estimates of Lemma 2.1.\fin

\noindent From now on any function belonging to  $H^1(\Omega)$ will be extended to a function belonging to 
$H^1(\widetilde{\Omega}_{\varepsilon,3})$. To make the notation simpler the extention of function  $\phi$  will still be  denoted 
$\phi$.
\medskip
\noindent In the sequel, we will make use of definitions and results from [2] and [5] concerning the periodic unfolding method.
Let us  recall the definition  of the unfolding operator ${\cal T}_\varepsilon$ which asociates a function ${\cal
T}_\varepsilon(\phi)\in L^1(\Omega\times Y)$ to each function  $\phi\in  L^1( \Omega_\varepsilon)$, 
$${\cal T}_\varepsilon(\phi)(x,y)=\phi\Bigl(\varepsilon\Bigr[{x\over \varepsilon}\Bigr]_Y+\varepsilon y\Bigr)\qquad \hbox{ for
$x\in \Omega$ and  $y\in Y$}.$$
\noindent We also recall  the approximate integration formula
$$\Bigl|\int_\Omega v-{1\over |Y|}\int_{\Omega\times Y}{\cal T}_\varepsilon(v)\Bigr|\le
||v||_{L^1(\widehat{\Omega}_{\varepsilon,1})}\qquad \forall v\in L^1(\Omega_\varepsilon)\leqno(2.2)$$ 
\noindent For the other properties of ${\cal T}_\varepsilon$, we refer the reader to [2] and [5].

\noindent Let $\phi\in H^1( \Omega)$ extended to $\widetilde{\Omega}_{\varepsilon,2}$. We have defined  the
scale-splitting operators ${\cal Q}_\varepsilon$ and ${\cal R}_\varepsilon$  (see [2]). The function ${\cal Q}_\varepsilon(\phi)$
is  the  restriction to $\Omega$ of  $Q_1$-interpolate of the discrete function $M^\varepsilon_Y(\phi)$ 
$$M^\varepsilon_Y(\phi)(x)={1\over |Y|}\int_Y \phi\Bigl(\varepsilon\Bigl[{x\over \varepsilon}\Bigr]+\varepsilon z\Bigr)
dz\qquad x\in \Omega$$  and ${\cal R}_\varepsilon(\phi)=\phi-{\cal Q}_\varepsilon(\phi)$. The operator ${\cal Q}_\varepsilon$
is linear and continuous from $H^1(\Omega)$ to $H^1(\Omega)$ and we have the estimates 
$$||{\cal Q}_\varepsilon(\phi)||_{H^1(\Omega)}\le C||\phi||_{H^1(\Omega)}\qquad ||\phi-{\cal Q}_\varepsilon(\phi)||_{
L^2(\Omega)}\le C\varepsilon||\nabla\phi||_{[L^2(\Omega)]^n}\qquad \forall \phi\in H^1(\Omega).$$ The constants depend 
on $n$ and $\partial\Omega$.
\medskip
\noindent{\bf Theorem 2.2 : }{\it Let $\phi$ be in $H^1(\Omega)$. There exists  $\widehat{\psi}_\varepsilon$ belonging to
$H^1_{per}(Y ; L^2(\Omega))$ such that
$$\left\{\eqalign{  &||\widehat{\psi}_\varepsilon||_{H^1(Y ; L^2(\Omega))}\le
C\bigl\{||\phi||_{L^2(\Omega)}+\varepsilon||\nabla \phi||_{[L^2(\Omega)]^n}\bigr\}\cr  &||{\cal
T}_\varepsilon(\phi)-\widehat{\psi}_\varepsilon||_{H^1(Y ; (H^1(\Omega ))^{'})}\le C\varepsilon\bigl\{||\phi||_{L^2(
\Omega)}+\varepsilon||\nabla \phi||_{[L^2(\Omega)]^n}\bigr\}\cr &\hskip 3.9cm+C\sqrt\varepsilon\bigl\{||\phi
||_{L^2(\widehat{\Omega}_{\varepsilon,2})}+
\varepsilon||\nabla \phi||_{[L^2(\widehat{\Omega}_{\varepsilon,2})]^n}\bigr\}\cr}\right.\leqno(2.3)$$
\noindent  The constants  depend only on $n$ and  $\partial\Omega$.}

\noindent{\bf Proof : } In this proof we use the same notation and the  same ideas as in Proposition 3.3 of [5].

\noindent Theorem 2.2 is proved in two steps. We reintroduce the unfolding operators ${\cal T}_{\varepsilon, i}$, defined in
[5], which for any $\phi\in H^1(\Omega)$, allow  us to estimate the difference between the restrictions to two
neighbouring cells of the unfolded of $\phi$ in $L^2(Y ; (H^1(\Omega))^{'})$. Then we evaluate the periodic defect of the
functions   $y\longrightarrow {\cal T}_\varepsilon(\phi)(.,y)$ thanks to Theorem 2.2 of [5].

\noindent Let  $\displaystyle K_i=Y\cup (\vec e_i+Y)$, $i\in\{1,\ldots,n\}$. For any $x$ in $\Omega$, the cell $\displaystyle 
\varepsilon\Bigl(\Bigl[{x\over \varepsilon}\Bigr]_Y+K_i\Bigr)$ is included in $\widetilde{\Omega}_{\varepsilon,2}$.

\noindent We recall that the unfolding  operator ${\cal T}_{\varepsilon, i}$ from $L^2(\widetilde{\Omega}_{\varepsilon,2})$
into $L^2(\Omega\times K_i)$ is defined  by
$$\forall\psi\in L^2(\widetilde{\Omega}_{\varepsilon,2}),\qquad {\cal T}_{\varepsilon,
i}(\psi)(x,y)=\psi\Bigl(\varepsilon\Bigl[{x\over
\varepsilon}\Bigr]_Y+ \varepsilon y\Bigr)\qquad \hbox{for $ x\in\Omega$ and  a. e. $ y\in  K_i$}.$$ 
\noindent  The restriction of  ${\cal T}_{ \varepsilon,i}(\psi)$ to $\Omega\times Y$ is equal to the unfolded  ${\cal
T}_\varepsilon(\psi)$ and we have the following equalities in $L^2(\Omega\times Y)$ :
$${\cal T}_{\varepsilon, i}(\psi)(.,..+\vec e_i)={\cal T}_{\varepsilon}(\psi)(.+\varepsilon \vec e_i,..),\qquad
i\in\{1,\ldots,n\}$$
\vskip 1mm
\noindent{\bf Step one.}  Let us take $\psi\in L^2(\widetilde{\Omega}_{\varepsilon,2})$. We evaluate the difference ${\cal
T}_{\varepsilon, i}(\psi)(.,.. +\vec e_i)-{\cal T}_{\varepsilon, i}(\psi)$ in $L^2(Y ; (H^1(\Omega))^{'})$.
\smallskip
\noindent For any $\Psi\in H^1(\Omega)$, extended  on $\widetilde{\Omega}_{\varepsilon, 1}$, a linear change of variables and
the relations above give
$$\eqalign{
\hbox{for a. e. } y\in Y,\qquad \int_{\Omega}{\cal T}_{\varepsilon, i}(\psi)(x,y+\vec e_i)\Psi(x)dx&=\int_{
\Omega}{\cal T}_{\varepsilon,i}(\psi)(x+\varepsilon\vec e_i,y)\Psi(x)dx\cr  &= \int_{ \Omega+\varepsilon\vec e_i}{\cal
T}_{\varepsilon,i}(\psi)(x,y)\Psi(x-\varepsilon\vec e_i)dx\cr}$$ We deduce 
$$\eqalign{   &\Bigl|\int_{\Omega}\bigl\{{\cal T}_{\varepsilon, i}(\psi)(., y+\vec e_i)-{\cal T}_{\varepsilon, i}(\psi)(.,
y)\bigr\}\Psi -\int_{\Omega}{\cal T}_{\varepsilon, i}(\psi)(., y)\bigl\{\Psi(.-\varepsilon\vec e_i)-\Psi\bigr\}
\Bigr|\cr
\le & C||{\cal T}_{\varepsilon,i}(\psi)(.,y)||_{L^2(\widehat{\Omega}_{\varepsilon,1})}||\Psi||_{
L^2(\widehat{\Omega}_{\varepsilon,1})}\qquad \hbox{for a. e. } y\in Y.\cr}$$ Since $\Omega$ is a bounded domain with
lipschitzian boundary and since $\Psi$ belongs to $ H^1(\widetilde{\Omega}_{\varepsilon,1})$ we  have
$$\left\{\eqalign{ &||\Psi||_{L^2(\widehat{\Omega}_{\varepsilon,1} )}\le C\sqrt\varepsilon\bigl\{||\Psi||_{L^2(\Omega)}+
||\nabla \Psi||_{[L^2(\Omega )]^n}\bigr\},\cr & ||\Psi(.-\varepsilon\vec e_i)-\Psi||_{L^2(\Omega)}\le
C\varepsilon\Bigl\|{\partial\Psi\over \partial x_i}\Bigr\|_{L^2(\Omega)},\qquad i\in\{1,\ldots,n\}, \cr}\right.\leqno(2.4)$$ hence
$$\eqalign{  <{\cal T}_{\varepsilon, i}&(\psi)(., y+\vec e_i)-{\cal T}_{\varepsilon, i}(\psi)(., y)\,,\,
\Psi>_{(H^1(\Omega))^{'},H^1( \Omega)}\cr  =&\int_{\Omega }\bigl\{{\cal T}_{\varepsilon, i}(\psi)(., y+\vec e_i)-{\cal
T}_{\varepsilon, i}(\psi)(., y)\bigr\}\Psi\cr
\le  &C\varepsilon||\nabla \Psi||_{[L^2(\Omega )]^n}||{\cal T}_{\varepsilon,i}(\psi)(.,y)||_{L^2(\Omega )}+  C\sqrt
\varepsilon||\Psi||_{H^1(\Omega)}||{\cal T}_{\varepsilon,i}(\psi)(.,y)||_{L^2(\widehat{
\Omega}_{\varepsilon,1})}.\cr}$$ We deduce
 that
$$||{\cal T}_{\varepsilon, i}(\psi)(.,y +\vec e_i)-{\cal T}_{\varepsilon, i}(\psi)(., y)||_{(H^1(\Omega))^{'}}\le  C
\varepsilon||{\cal T}_{\varepsilon,i}(\psi)(.,y)||_{L^2( \Omega)}+  C\sqrt \varepsilon||{\cal T}_{\varepsilon,i}(\psi)(.,y)
||_{L^2(\widehat{\Omega}_{\varepsilon,1})}.$$
\noindent Which leads to the following estimate of the difference between ${\cal T}_{\varepsilon, i}(\psi)_{|_{\Omega\times
Y}}$ and one of its translated :
$$||{\cal T}_{\varepsilon, i}(\psi)(.,.. +\vec e_i)-{\cal T}_{\varepsilon, i}(\psi)||_{L^2(Y ; (H^1(\Omega))^{'})}
\le   C\varepsilon||\psi||_{L^2(\widetilde{\Omega}_{\varepsilon,3})}+  C\sqrt
\varepsilon||\psi||_{L^2(\widehat{\Omega}_{\varepsilon,2})}.\leqno(2.5)$$ The constant depends only on $n$ and on the
boundary of  $\Omega$.
\medskip
\noindent{\bf Step two. } Let $\phi\in H^1(\Omega)$. The estimate $(2.5)$ applied to $\phi$ and its partial derivatives give  us
$$\eqalign{  ||{\cal T}_{\varepsilon, i}(\phi)(., .. +\vec e_i)-{\cal T}_{\varepsilon, i}(\phi) ||_{ L^2(Y ; (H^1(\Omega))^{'})}
& \le C\varepsilon\bigl\{||\phi||_{L^2(\Omega)}+\varepsilon ||\nabla \phi||_{[L^2(\Omega)]^n}\bigr\}  +C\sqrt \varepsilon||\phi
||_{L^2(\widehat{\Omega}_{\varepsilon,2})}\cr ||{\cal T}_{\varepsilon, i}(\nabla \phi)(., .. +\vec e_i)-{\cal T}_{\varepsilon,
i}(\nabla  \phi)||_{ [L^2(Y ; (H^1 (\Omega))^{'}]^n)} &\le C\bigl\{\varepsilon||\nabla \phi||_{[L^2(\Omega)]^n}+\sqrt
\varepsilon||\nabla \phi||_{[L^2(\widehat{\Omega}_{\varepsilon,2})]^n}\bigr\}\cr}$$ We recall  that 
$\nabla_y\bigl({\cal T}_{\varepsilon, i}(\phi)\bigr)=\varepsilon{\cal T}_{\varepsilon, i}(\nabla \phi)$ (see [3]). The above
estimates can also be written as follows :
$$\eqalign{ ||{\cal T}_{\varepsilon,i}(\phi)(., .. +\vec e_i)-{\cal T}_{\varepsilon,i}(\phi)||_{H^1(Y ; (H^1(\Omega))^{'})} &\le
C\varepsilon\bigl\{||\phi||_{L^2(\Omega)}+\varepsilon||\nabla \phi||_{[L^2(\Omega)]^n}+\sqrt
\varepsilon||\nabla \phi||_{[L^2(\widehat{\Omega}_{\varepsilon,2})]^n}\bigr\}\cr &+C\sqrt \varepsilon||\phi
||_{L^2(\widehat{\Omega}_{\varepsilon,2})}\cr}$$ 
\noindent From these inequalities, for any
$i\in\{1,\ldots,n\}$,  we deduce the estimate of the difference of the traces of
 $y\longrightarrow {\cal T}_{\varepsilon}(\phi)(.,y)$ on the faces $Y_i$ and $\vec e_i+Y_i$ 
$$\left\{\eqalign{ ||{\cal T}_{\varepsilon}(\phi)(.,.. +\vec e_i)-{\cal T}_{\varepsilon}(\phi)||_{H^{1/ 2}(Y_i  ;
(H^1(\Omega))^{'})}&\le C\varepsilon\bigl\{||\phi||_{L^2(\Omega)}+\varepsilon||\nabla \phi||_{[L^2(\Omega)]^n}\bigr\}\cr
&+C\sqrt \varepsilon\bigl\{||\phi ||_{L^2(\widehat{\Omega}_{\varepsilon,2})}+
\varepsilon||\nabla \phi||_{[L^2(\widehat{\Omega}_{\varepsilon,2})]^n}\bigr\}\cr}\right. $$ It measures the periodic defect
of $y\longrightarrow{\cal T}_\varepsilon(\phi)(.,y)$. We decompose  ${\cal T}_{\varepsilon}(\phi)$ into the sum of an element
belonging to $H^1_{per}(Y ; L^2(\Omega ))$ and an element belonging to $\bigl(H^1(Y ; L^2(\Omega))\bigr)^{\perp}$ (the
orthogonal of $H^1_{per}(Y ; L^2(\Omega))$ in $H^1(Y ; L^2(\Omega))$, see [5])
$${\cal T}_{\varepsilon}(\phi)=\widehat{\psi}_\varepsilon+\overline{\phi}_\varepsilon,\qquad\widehat{\psi}_\varepsilon
\in H^1_{per}(Y ; L^2(\Omega )),\qquad \overline{\phi}_\varepsilon\in  \bigl(H^1(Y ; L^2(\Omega))\bigr)^\perp\leqno(2.6)$$
From the Riesz Theorem the dual space $(H^1(\Omega))^{'}$ is a Hilbert space isomorphic to $H^1(\Omega)$. The function
$y\longrightarrow {\cal T}_\varepsilon(
\phi)(., y)$ takes its values in a finite dimensionnal  space, 
$$\overline{\phi}_\varepsilon(., ..) =\sum_{\xi\in \Xi_\varepsilon}\overline{\phi}_{\varepsilon, \xi}(..)\chi_\xi(.)$$ where
$\chi_\xi(.)$ is the  characteristic  function of the cell $\varepsilon(\xi+Y)$ and where $\overline{\phi}_{\varepsilon, \xi}(..)\in
\bigl(H^1(Y)\bigr)^\perp$  (the orthogonal of $H^1_{per}(Y)$ in $H^1(Y)$, see [5]). Hence the   decomposing $(2.6)$ is the
same in $H^1(Y ; (H^1(\Omega))^{'})$.  As the decomposing is orthogonal, we have
$$\eqalign{ ||\widehat{\psi}_\varepsilon||^2_{H^1(Y ; L^2(\Omega))}+||\overline{\phi}_\varepsilon||^2_{H^1(Y ; L^2(\Omega))}
=  ||{\cal T}_\varepsilon(\phi)||^2_{H^1(Y ; L^2(\Omega))} 
\le  C\bigl\{||\phi||_{L^2(\Omega)}+\varepsilon||\nabla \phi||_{[L^2(\Omega)]^n}\bigr\}^2\cr}$$ Hence we have the first
inequality
$(2.3)$ and an estimate of $\overline{\phi}_\varepsilon$ in
$H^1(Y ; L^2(\Omega))$. From Theorem 2.2 of [5] and  $(2.5)$ we obtain a finer estimate of $\overline{\phi}_\varepsilon$ in
$H^1(Y ; (H^1(\Omega))^{'})$
$$||\overline{\phi}_\varepsilon||_{H^1(Y ; (H^1(\Omega))^{'})} \le 
C\varepsilon\bigl\{||\phi||_{L^2(\Omega)}+\varepsilon||\nabla \phi||_{[L^2(\Omega)]^n}+\sqrt
\varepsilon||\nabla \phi||_{[L^2(\widehat{\Omega}_{\varepsilon,2})]^n}\bigr\}+C\sqrt \varepsilon||\phi
||_{L^2(\widehat{\Omega}_{\varepsilon,2})}$$ It is the second inequality in $(2.3)$.\fin
\noindent{\bf Theorem 2.3 : }{\it For any  $\phi\in H^1(\Omega)$, there exists
$\widehat{\phi}_\varepsilon\in H^1_{per}(Y ; L^2(\Omega))$ such that
$$\left\{\eqalign{&||\widehat{\phi}_\varepsilon||_{H^1(Y ; L^2(\Omega))}\le C||\nabla \phi||_{[L^2(\Omega)]^n},\cr  &||{\cal
T}_\varepsilon(\nabla_x \phi)-\nabla_x \phi-\nabla_y\widehat{\phi}_\varepsilon||_{ [L^2(Y ; (H^1(\Omega))^{'})]^n}
\le C\varepsilon||\nabla \phi||_{[L^2(\Omega)]^n}+C\sqrt\varepsilon||\nabla \phi||_{[L^2(\widehat{\Omega}_{\varepsilon,3})
]^n}.}\right.\leqno(2.7)$$ 
\noindent  The constants  depend only on $n$ and  $\partial\Omega$.}

\noindent{\bf Proof : } Let $\phi\in H^1(\Omega)$. The function $\phi$ is decomposed 
$$\phi=\Phi+\varepsilon \underline{\phi} ,
\quad\hbox{where}\enskip \Phi={\cal Q}_\varepsilon(\phi)\quad\hbox{and} \quad
\underline{\phi}={1\over
\varepsilon}{\cal R}_\varepsilon(\phi).$$ with the following estimate :
$$||\nabla \Phi||_{[L^2(\Omega)]^n}+||\underline{\phi}||_{L^2(\Omega)}+\varepsilon||\nabla 
\underline{\phi}||_{[L^2(\Omega)]^n}\le C||\nabla \phi||_{[L^2(\Omega)]^n}.\leqno(2.8)$$ We apply the PoincarŽ-Wirtinger
inequality to the function $\phi$ in each cell of the form $\varepsilon(\xi +K_i)$ and of the form $\varepsilon(\xi+Y)$ included
in   $\widehat{\Omega}_{\varepsilon,3}$. We deduce that
$$\eqalign{ &||\nabla {\cal Q}_\varepsilon(\phi)||_{[L^2(\widehat{\Omega}_{\varepsilon,2})]^n}\le C||\nabla
\phi||_{[L^2(\widehat{\Omega}_{\varepsilon,3})]^n}\cr
\Longrightarrow \quad& ||\nabla \underline{\phi}||_{[L^2(\widehat{\Omega}_{\varepsilon,2})]^n}\le {C\over 
\varepsilon}||\nabla \phi||_{[L^2(\widehat{\Omega}_{\varepsilon,3})]^n}\cr}$$ We also have (see [3])
$$||\underline{\phi}||_{L^2(\widehat{\Omega}_{\varepsilon,2})}={1\over \varepsilon}||\phi-{\cal
Q}_\varepsilon(\phi)||_{L^2(\widehat{\Omega}_{\varepsilon,2})}\le C||\nabla
\phi||_{[L^2(\widehat{\Omega}_{\varepsilon,3})]^n}$$

\noindent Theorem 3 applied to $\underline{\phi}$ gives us the existence of an element $\widehat{\phi}_{\varepsilon}$ in
$H^1_{per}(Y ; L^2(\Omega))$ such that
$$\left\{\eqalign{  & ||\widehat{\phi}_{\varepsilon}||_{H^1(Y ; L^2(\Omega))}\le C||\nabla \phi||_{[L^2(\Omega)]^n},\cr
 & ||{\cal T}_\varepsilon(\underline{\phi})-\widehat{\phi}_{\varepsilon}||_{H^1(Y ;  (H^1(\Omega))^{'})}\le C\varepsilon||\nabla
\phi||_{[L^2(\Omega)]^n}+ {C\sqrt\varepsilon}||\nabla
\phi||_{[L^2(\widehat{\Omega}_{\varepsilon,3})]^n}.\cr}\right.\leqno(2.9)$$

\noindent   We  evaluate $||{\cal T}_\varepsilon(\nabla \Phi)-\nabla \Phi||_{[L^2(Y ; (H^1(\Omega))^{'})]^n}$. 
\vskip 1mm
\noindent From Lemma 2.2 we have
$$\Bigl\|{\partial\Phi\over \partial x_i}-M^\varepsilon_Y\Bigl({\partial\Phi\over \partial x_i}\Bigr)\Bigr\|_{(H^1(\Omega))^{'}}\le
 C\varepsilon||\nabla \phi||_{[L^2(\Omega)]^n}+ C\sqrt
\varepsilon||\nabla\phi||_{L^2(\widehat{\Omega}_{\varepsilon,3})]^n}\leqno(2.10)$$ From the definition of
$\Phi$ it results that $\displaystyle y\longrightarrow{\cal T}_\varepsilon\Bigl({\partial\Phi\over
\partial x_i}\Bigr)(.,y)$ is linear with respect to each variable.  For any $\psi\in H^1(\Omega)$, we have
$$\eqalign{ <{\cal T}_\varepsilon \Bigl({\partial\Phi\over \partial x_1}\Bigr)(., y)-M^\varepsilon_Y\Bigl({\partial\Phi\over\partial
x_1}\Bigr),\psi>_{(H^1(\Omega))^{'}, H^1(\Omega)}&=
\int_\Omega\Bigl\{{\cal T}_\varepsilon\Bigl({\partial\Phi\over \partial x_1}\Bigr)(., y)-M^\varepsilon_Y\Bigl({\partial\Phi\over
\partial x_1}\Bigr)\Bigr\}\psi\cr &=\int_{\Omega_\varepsilon}\Bigl\{{\cal T}_\varepsilon\Bigl({\partial\Phi\over \partial
x_1}\Bigr)(., y)-M^\varepsilon_Y\Bigl({\partial\Phi\over
\partial x_1}\Bigr)\Bigr\}M^\varepsilon_Y(\psi)\cr
 &+\int_{\Omega\setminus\Omega_\varepsilon}\Bigl\{{\cal T}_\varepsilon\Bigl({\partial\Phi\over \partial x_1}\Bigr)(.,
y)-M^\varepsilon_Y\Bigl({\partial\Phi\over\partial x_1}\Bigr)\Bigr\} \psi \cr}$$ We have
$$\int_{\Omega\setminus\Omega_\varepsilon}\Bigl\{{\cal T}_\varepsilon\Bigl({\partial\Phi\over \partial x_1}\Bigr)(., y)
-M^\varepsilon_Y\Bigl({\partial\Phi\over\partial x_1}\Bigr)\Bigr\} \psi \le C\sqrt \varepsilon||\nabla\phi||_{L^2(\widehat{
\Omega}_{\varepsilon,3})]^n}\bigl\{||\psi||_{L^2(\Omega)}+||\nabla\psi||_{L^2( \Omega)]^n}\bigr\}$$ Besides, as in
Theorem 3.4 of [5] we show  that
$$\eqalign{
\int_{\Omega_\varepsilon}\Bigl\{{\cal T}_\varepsilon\Bigl({\partial\Phi\over \partial x_1}\Bigr)(., y)-M^\varepsilon_Y
\Bigl({\partial\Phi\over\partial x_1}\Bigr)\Bigr\}M^\varepsilon_Y(\psi)&\le  C\varepsilon||\nabla \phi||_{[L^2(
\Omega)]^n}||\nabla\psi||_{[L^2(\Omega)]^n}\cr &+
C\sqrt\varepsilon||\nabla\phi||_{L^2(\widehat{\Omega}_{\varepsilon, 3})]^n}\bigl\{||\psi||_{L^2(\Omega)}+||\nabla\psi||_{
L^2(\Omega)]^n}\bigr\}\cr}$$ and eventually
$$\forall y\in Y,\qquad \Bigl\|{\cal T}_\varepsilon\Bigl({\partial \Phi\over \partial x_1}\Bigr)(., y)- M^\varepsilon_Y\Bigl({\partial
\Phi\over \partial x_1}\Bigr) \Bigr\|_{(H^1(\Omega))^{'}}\le C\varepsilon ||\nabla \phi||_{[L^2(\Omega)]^n}+ C\sqrt
\varepsilon||\nabla\phi||_{[L^2(\widehat{\Omega}_{\varepsilon,3})]^n}.$$ Considering $(2.10)$ and all the partial derivatives,
we obtain
$$||{\cal T}_\varepsilon(\nabla \Phi)-\nabla \Phi||_{[L^2(Y ; (H^1(\Omega))^{'})]^n}\le C\varepsilon||\nabla \phi||_{
[L^2(\Omega)]^n}+C\sqrt \varepsilon||\nabla\phi||_{[L^2(\widehat{\Omega}_{\varepsilon,3})]^n}$$ Moreover we have
$$\eqalign{
\int_\Omega {\partial \underline{\phi}\over \partial x_i}\psi=\int_{\partial\Omega}\underline{\phi}n_i\psi
-\int_\Omega\underline{\phi}{\partial \psi\over \partial x_i}\le
C\bigl\{||\underline{\phi}||_{L^2(\partial\Omega)}+C||\underline{\phi}||_{L^2(\Omega)}\bigr\}||\psi||_{H^1(\Omega)}\cr
||\underline{\phi}||_{L^2(\partial\Omega)}\le {C\over \sqrt\varepsilon} ||
\underline{\phi}||_{L^2(\widehat{\Omega}_{\varepsilon,1})]}+C\sqrt\varepsilon ||\nabla
\underline{\phi}||_{[L^2(\widehat{\Omega}_{\varepsilon,1})]^n}\le {C\over \sqrt \varepsilon} ||\nabla
\phi||_{[L^2(\widehat{\Omega}_{\varepsilon,3})]^n}\cr}$$ hence
$||\varepsilon\nabla \underline{\phi}||_{[(H^1(\Omega))^{'}]^n}\le  C\varepsilon||\nabla \phi||_{[L^2(\Omega)]^n}+C\sqrt
\varepsilon ||\nabla\phi||_{[L^2(\widehat{\Omega}_{\varepsilon,3})]^n}$. Thanks to $(2.9)$ and  to the
above inequalities  the second estimate of $(2.7)$ is proved.\fin
\medskip
\noindent{\bf 2.2 Projection theorems in $L^2(Y; (H^s(\Omega))^{'})$, $0<s<1$.}
\medskip
\noindent The space $H^s(\Omega)$, $0<s<1$, is defined by
$$H^s(\Omega)=\Bigl\{\phi\in L^2(\Omega)\;|\; \int_{\Omega\times \Omega}{|\phi(x)-\phi(x{'})|^2\over
|x-x^{'}|^{n+2s}}dxdx^{'}<+\infty\Bigr\}.$$ Equipped with the inner product
$$<\phi,\psi>_s=\int_\Omega\phi\,\psi+\int_{\Omega\times
\Omega}{\bigl(\phi(x)-\phi(x{'})\bigr)\bigl(\psi(x)-\psi(x^{'})\bigr)
\over |x-x^{'}|^{n+2s}}dxdx^{'}$$  $H^s(\Omega)$ is a Hilbert separable space. We denote $||.||_{s,\Omega}$ the norm
associated to this inner product.

\noindent As we have done in Lemma 2.1 we build a linear and continuous extension operator 
${\cal P}$ from $H^s(\Omega)$, $0<s<1$, into  $H^s(\widetilde{\Omega}_{\varepsilon, 4})$ verifying
$$||{\cal P}(\phi)||_{s,\widetilde{\Omega}_{\varepsilon, 4}}\le C || \phi||_{s,\Omega}$$ \noindent  The constant  depends only on
$n$, $s$ and  $\partial\Omega$.

\noindent From now on any function belonging to  $H^s(\Omega)$ will be extended to a function belonging to 
$H^s(\widetilde{\Omega}_{\varepsilon,3})$, $0<s<1$. To make the notation simpler the extention of function  $\phi$  will still
be  denoted  $\phi$.
\medskip
\noindent{\bf Lemma 2.4 : }{\it  For any $\phi\in H^s(\Omega)$, $0<s<1$, we have
$$\left\{\eqalign{
&||\phi-M^\varepsilon_Y(\phi)||_{L^2( \Omega)}\le C\varepsilon^s||\phi||_{s,\Omega}\qquad ||\phi-{\cal
Q}_\varepsilon(\phi)||_{L^2( \Omega)}\le C\varepsilon^s||\phi||_{s,\Omega}\cr 
&||\nabla{\cal Q}_\varepsilon(\phi)||_{[L^2(\widetilde{\Omega}_{\varepsilon,3})]^n}\le C\varepsilon^{s-1}||\phi||_{s,
\Omega},\quad ||\phi||_{L^2(\widehat{\Omega}_{\varepsilon,3})}\le C\varepsilon^{s/2}||\phi||_{s,\Omega}\cr
&||\phi-\phi(.+\varepsilon\vec e_i)||_{L^2( \Omega)}\le C\varepsilon^s||\phi||_{s,\Omega},\qquad i\in\{1,\ldots,n\}\cr 
& ||{\cal Q}_\varepsilon(\phi)||_{L^2(\partial\Omega)}\le C\varepsilon^{(s-1)/2}||\phi||_{s,\Omega}
\cr}\right.\leqno(2.11)$$ The constants depend on $n$, $s$ and  $\partial\Omega$.}

\noindent{\bf Proof : } For any $\psi$ belonging to  $H^s(Y)$, $0<s<1$,  we have the Poincar\'e-Wirtinger inequality
$$||\psi-M_Y(\psi)||_{L^2(Y)}\le C||\psi||_{s,Y}$$ where $M_Y(\psi)$ is the mean of $\psi$ in the cell $Y$. The constant
depends only on  $n$. We  immediately deduce the upper
bound $||\phi-M^\varepsilon_Y(\phi)||_{L^2(\widetilde{\Omega}_{\varepsilon,4})}\le C\varepsilon^s||\phi||_{s,\Omega}$. 
We apply the PoincarŽ-Wirtinger inequality  to the restriction of $\phi$  to  two neighbouring cells included in 
$\widetilde{\Omega}_{\varepsilon,4}$ and we obtain the estimate of the  gradient of
${\cal Q}_\varepsilon(\phi)$  in $\widetilde{\Omega}_{\varepsilon,3}$ ($||\nabla{\cal
Q}_\varepsilon(\phi)||_{[L^2(\widetilde{\Omega}_{\varepsilon,3})]^n}\le C\varepsilon^{s-1}||\phi||_{s,\Omega}$)  and then the
upper bound   $||\phi-{\cal Q}_\varepsilon(\phi)||_{L^2(\widetilde{\Omega}_{\varepsilon,3})}\le
C\varepsilon^s||\phi||_{s,\Omega}$ thanks to  the estimate of
$||\phi-M^\varepsilon_Y(\phi)||_{L^2(\widetilde{\Omega}_{\varepsilon,3})}$.  The function ${\cal Q}_\varepsilon(\phi)$ belongs
to $H^1(\widetilde{\Omega}_{\varepsilon,3})$, hence considering a neighbourhood of
$\partial\widetilde{\Omega}_{\varepsilon,3}$  (included in $\widetilde{\Omega}_{\varepsilon,3}$) of thickness
$\varepsilon^{1-s}$ we show that 
$$||{\cal Q}_\varepsilon(\phi)||_{L^2(\widehat{\Omega}_{\varepsilon, 3})}\le  C\varepsilon^{s/2}||\phi||_{s,\Omega}\qquad
\Longrightarrow \qquad ||\phi||_{L^2(\widehat{\Omega}_{\varepsilon, 3})}\le 
C\varepsilon^{s/2}||\phi||_{s,\Omega}.$$ We have
$$\eqalign{ ||\phi-\phi(.+\varepsilon\vec e_i)||_{L^2( \Omega)} &\le ||\phi-{\cal Q}_\varepsilon(\phi)||_{L^2( \Omega)}+||{\cal
Q}_\varepsilon(\phi)-{\cal Q}_\varepsilon(\phi)(.+\varepsilon\vec e_i)||_{L^2( \Omega)}\cr &+||\phi(.+\varepsilon\vec e_i)-{\cal
Q}_\varepsilon(\phi)(.+\varepsilon\vec e_i)||_{L^2(\Omega)}
\le C\varepsilon^s||\phi||_{s,\Omega}\cr}$$ thanks to the upper bounds of
$||\phi-{\cal Q}_\varepsilon(\phi)||_{L^2(\widetilde{\Omega}_{\varepsilon,3})}$ and $||\nabla{\cal
Q}_\varepsilon(\phi)||_{[L^2(\widetilde{\Omega}_{\varepsilon,3})]^n}$. The last inequality of the lemma is the
consequence  of the estimates $||\nabla{\cal Q}_\varepsilon(\phi)||_{[L^2(\widetilde{\Omega}_{\varepsilon,3})]^n}\le
C\varepsilon^s||\phi||_{s,\Omega}$ and $||{\cal Q}_\varepsilon(\phi)||_{ L^2(\widetilde{\Omega}_{\varepsilon,3})}\le
C||\phi||_{s,\Omega}$.\fin
\noindent{\bf Corollary : }{\it For any $s\in ]0,1[$ and for any  $\phi\in H^s(\Omega)$ we have
$$\left\{\eqalign{ ||{\cal Q}_\varepsilon(\phi)-M^\varepsilon_Y(\phi)||_{L^2(\Omega)}\le C\varepsilon^s||\phi||_{s,\Omega}\cr
||\phi-{\cal T}_\varepsilon(\phi)||_{L^2(\Omega\times Y)}\le C\varepsilon^s||\phi||_{s,\Omega}\cr}\right.\leqno(2.12) $$ The
constants depend on $n$, $s$ and  $\partial\Omega$.}

\noindent{\bf Proof : } The inequalities $(2.12) $ are the consequences of $(2.11)$.\fin

\noindent{\bf Theorem 2.5 : }{\it Let $\phi$ be in $H^1(\Omega)$. There exists  $\widehat{\psi}_\varepsilon$ belonging to
$H^1_{per}(Y ; L^2(\Omega))$ such that for any $s\in ]0,1[$
$$\left\{\eqalign{  &||\widehat{\psi}_\varepsilon||_{H^1(Y ; L^2(\Omega))}\le
C\bigl\{||\phi||_{L^2(\Omega)}+\varepsilon||\nabla \phi||_{[L^2(\Omega)]^n}\bigr\}\cr   &||{\cal
T}_\varepsilon(\phi)-\widehat{\psi}_\varepsilon||_{H^1(Y ; (H^s(\Omega ))^{'})}\le
C\varepsilon^s\bigl\{||\phi||_{L^2(\Omega)}+\varepsilon||\nabla \phi||_{[L^2(\Omega)]^n}\bigr\}\cr &\hskip
3.9cm+C\varepsilon^{s/2}\bigl\{||\phi ||_{L^2(\widehat{\Omega}_{\varepsilon,2})}+
\varepsilon||\nabla \phi||_{[L^2(\widehat{\Omega}_{\varepsilon,2})]^n}\bigr\}\cr}\right.\leqno(2.13) $$
\noindent  The constants  depend only on $n$, $s$ and  $\partial\Omega$.}

\noindent{\bf Proof : } With a few modifications we prove Theorem 2.5 as Theorem 2.2. Thanks to Lemma 2.4 we replace the
inequalities $(2.4)$ of step one in Theorem 2.2 by
$$\forall\Psi\in H^s(\Omega),\qquad\left\{\eqalign{ &||\Psi||_{L^2(\widehat{\Omega}_{\varepsilon,1} )}\le C
\varepsilon^{s/2}||\Psi||_{s,\Omega},\cr & ||\Psi(.-\varepsilon\vec e_i)-\Psi||_{L^2(\Omega)}\le
C\varepsilon^s||\Psi||_{s,\Omega},\qquad i\in\{1,\ldots,n\}.\cr}\right.$$\fin
\noindent{\bf Theorem 2.6 : }{\it For any  $\phi\in H^1(\Omega)$, there exists
$\widehat{\phi}_\varepsilon\in H^1_{per}(Y ; L^2(\Omega))$ such that
$$\left\{\eqalign{&||\widehat{\phi}_\varepsilon||_{H^1(Y ; L^2(\Omega))}\le C||\nabla \phi||_{[L^2(\Omega)]^n},\cr  &||{\cal
T}_\varepsilon(\nabla_x \phi)-\nabla_x \phi-\nabla_y\widehat{\phi}_\varepsilon||_{ [L^2(Y ; (H^s(\Omega) )^{'})]^n}\le 
C\varepsilon^s||\nabla \phi||_{[L^2(\Omega)]^n}+C\varepsilon^{s/2}||\nabla \phi||_{
[L^2(\widehat{\Omega}_{\varepsilon,3})]^n}.}\right.\leqno(2.14)$$ 
\noindent  The constants  depend only on $n$, $s$ and  $\partial\Omega$.}

\noindent{\bf Proof : }  With a few modifications we prove Theorem 2.6 as Theorem 2.3. Proceeding as Theorem
3.4  in [5] and thanks to Lemma 2.4, we show that 
$$ ||{\cal T}_\varepsilon(\nabla \Phi)-\nabla \Phi||_{[L^2(Y ; (H^s(\Omega))^{'})]^n}\le C\varepsilon^s ||\nabla
\phi||_{[L^2(\Omega)]^n}+ C \varepsilon^{s/2}||\nabla\phi||_{[L^2(\widehat{\Omega}_{\varepsilon,3})]^n}$$  where
$\phi=\Phi+\varepsilon\underline{\phi}$, $\Phi={\cal Q}_\varepsilon(\phi)$. Now let $\psi$ be in $H^s(\Omega)$. We have
$$\eqalign{
\int_\Omega{\partial \underline{\phi}\over \partial x_i}\psi &=\int_\Omega{\partial \underline{\phi}\over \partial x_i}\bigl(
\psi-{\cal Q}_\varepsilon(\psi)\bigr)+\int_\Omega{\partial \underline{\phi}\over \partial x_i}{\cal Q}_\varepsilon(\psi)=
\int_\Omega{\partial \underline{\phi}\over \partial x_i}\bigl(\psi-{\cal
Q}_\varepsilon(\psi)\bigr)+\int_{\partial\Omega} \underline{\phi}n_i{\cal Q}_\varepsilon(\psi)-\int_\Omega
\underline{\phi}{\partial {\cal Q}_\varepsilon(\psi)\over \partial x_i}\cr
&\le ||\nabla \underline{\phi}||_{[L^2(\Omega)]^n}||\psi-{\cal Q}_\varepsilon(\psi)||_{L^2(\Omega)}+
||\underline{\phi}||_{L^2(\partial\Omega)}||{\cal Q}_\varepsilon(\psi)||_{L^2(\partial\Omega)}
+||\underline{\phi}||_{L^2(\Omega)}||{\cal Q}_\varepsilon(\psi)||_{[L^2(\Omega)]^n}\cr}$$
 hence $||\varepsilon\nabla \underline{\phi}||_{[(H^s(\Omega))^{'}]^n}\le C\varepsilon^s ||\nabla \phi||_{[L^2(\Omega)]^n}+
C \varepsilon^{s/2}||\nabla\phi||_{[L^2(\widehat{\Omega}_{\varepsilon,3})]^n}$ thanks to the estimates of
$\underline{\phi}$ (see Theorem 2.3) and the  inequalities of Lemma 2.4. \fin
\vfill\eject
\noindent{\bf 3. Error estimate in the classical homogenization problem}
\medskip
\noindent We consider the following homogenization problem :
$$ \left\{\eqalign{
& \phi^{\varepsilon}\in H^1_{\Gamma_0}(\Omega),\cr
&\int_\Omega A\bigl(\bigl\{{.\over\varepsilon}\bigr\}\bigr)\nabla\phi^{\varepsilon}.\nabla u=
\int_{\Omega}f u,\cr
&\forall u\in H^1_{\Gamma_0}(\Omega),\cr}\right.\leqno(3.1)$$  where  

$\bullet$ $\Omega$ is a bounded domain in $\R^n$ with lipschitzian boundary,

$\bullet$  $\Gamma_0$ is a measurable set of $\partial \Omega$ with measure nonnull or  $\Gamma_0=\emptyset$, 

$\bullet$ $H^1_{\Gamma_0}(\Omega)=\bigl\{\phi\in H^1(\Omega)\; |\;
\phi=0\; \hbox{on}\;\Gamma_0\bigr\}$,

$\bullet$ $f\in L^2(\Omega)$,

$\bullet$ $A$ is a square matrix of elements belonging to $L^\infty(Y)$, verifying the condition of
uniform ellipticity  $c|\xi|^2\le A(y)\xi.\xi\le C|\xi|^2$  a.e. $y\in Y$, with $c$ and  $C$ strictly positive constants.  
\medskip
\noindent If $\Gamma_0=\emptyset$, we suppose that $\displaystyle \int_\Omega f=\int_\Omega \phi^\varepsilon=0$
\medskip
We have shown, see [2], that $\nabla \phi^\varepsilon-\nabla\Phi-{\cal U}_\varepsilon\bigl(\nabla_y
\widehat{\phi}\bigr)$  strongly converges towards $0$ in $[L^2(\Omega)]^n$, where ${\cal U}_\varepsilon$ is the averaging
operator  defined by
$$V\in L^2(\Omega\times Y)\qquad {\cal U}_\varepsilon(V)(x)=\int_YV\Bigl(\varepsilon\Bigr[{x\over
\varepsilon}\Bigr]+\varepsilon z, \Bigl\{{x\over
\varepsilon}\Bigr\}\Bigr)dz,\qquad {\cal U}_\varepsilon(V)\in L^2(\Omega),$$ and where  
$$(\Phi,\widehat{\phi})\in  H^1_{\Gamma_0}(\Omega)\times L^{2}(\Omega,H^1_{per}(Y)/\R)$$ is the solution of the
limit problem of unfolding homogenization
 $$\left\{\eqalign{
 &\forall (U,\widehat{u})\in H^1_{\Gamma_0}(\Omega)\times L^{2}(\Omega ; H^1_{per}(Y)/\R)\cr&
\int_{\Omega}\int_{Y}A\bigl\{\nabla_x \Phi+\nabla_y\widehat{\phi}\bigr\}.
\bigl\{\nabla_x U+\nabla_y\widehat{u}\bigr\}=\int_\Omega f U.}\right.\leqno(3.2)$$ 
\noindent If $\Gamma_0=\emptyset$, we take $\displaystyle \int_\Omega\Phi=0$.

\noindent   We recall that the correctors $\chi_i$, $i\in\{1,\ldots,n\}$, are the solutions of the following variational problems :
\vskip-0.2cm
$$\chi_i\in H^1_{per}(Y),\qquad\int_Y\chi_i=0,\qquad  \int_YA(y)\nabla_y\bigl(\chi_i(y)+y_i\bigr)\nabla_y\psi(y)dy=0,\qquad
\forall\psi\in H^1_{per}(Y)$$
\noindent They allow us to express $\widehat{\phi}$ in terms of $\nabla \Phi$
$$\widehat{\phi}=\sum_{i=1}^n{\partial\Phi\over\partial x_i}\chi_i,$$ and to give the  homogenized problem verified by
$\Phi$ 
$$\int_\Omega{\cal A}\nabla \Phi\nabla U=\int_\Omega f U\qquad\qquad \forall U\in
H^1_{\Gamma_0}(\Omega) \leqno(3.3)$$ where (see [3])
$${\cal A}_{ij}={1\over |Y|}\sum_{k,l=1}^n\int_Y a_{kl}{\partial(y_i+\chi_i)\over \partial y_k}{\partial(y_j+\chi_j)\over
\partial y_l}.$$
\noindent{\bf 3.1 First case : smooth boundary and homogeneous Dirichlet or Neumann limits conditions}
\medskip
\noindent{\it  In this paragraph we suppose  that

$\bullet$ $\Omega$ is a bounded domain in $\R^n$ with ${\cal C}^{1,1}$ boundary,

$\bullet$  $\Gamma_0=\partial\Omega$ (homogeneous Dirichlet  condition) or  $\Gamma_0=\emptyset$ (homogeneous
Neumann condition). }
\smallskip
\noindent In Theorems 4.1  and 4.2 in [5]   we  gave  the following error estimate for the solution of problem $(3.1)$ :
$$||\phi^\varepsilon-\Phi||_{L^2(\Omega)}+||\nabla \phi^\varepsilon-\nabla \Phi- \sum_{i=1}^n
 {\cal Q}_\varepsilon\Bigl({\partial\Phi\over\partial x_i}\Bigr)\nabla_y\chi_i\bigl(\bigl\{{.\over \varepsilon}\bigr\}\bigr)||_{[L^2(\Omega)]^n}\le C\varepsilon^{ 1/2}||f||_{L^2(\Omega)},\leqno(3.4)$$ the constant depends on $n$,
$A$  and $\partial\Omega$. In Theorem 3.2 we are going to complete these estimates.
\medskip
\noindent{\bf Lemma 3.1 : }{\it  We have
$$  ||\nabla \phi^\varepsilon||_{[L^2(\widehat{\Omega}_{\varepsilon,3})]^n}\le
C\sqrt\varepsilon||f||_{L^2(\Omega)}\leqno(3.5)$$ The constant depends on
$n$, $A$ and $\partial\Omega$.}

\noindent{\bf Proof : } The boundary of $\Omega$ being of class ${\cal C}^{1,1}$ we deduce that the   solution $\Phi$ of the
homogenized problem $(3.3.i)$ belongs to $H^2(\Omega)$ and verifies $||\Phi||_{H^2(\Omega)}\le C||f||_{L^2(\Omega)}$. The
estimate of Lemma 3.1 is a consequence of  $(2.1)$, and of $(3.4)$ and of the following inequality :
$$\eqalign{ 
|| \nabla \Phi- \sum_{i=1}^n {\cal Q}_\varepsilon\Bigl({\partial\Phi\over\partial x_i}\Bigr)\nabla_y\chi_i\bigl( {.\over
\varepsilon}\bigr)||_{[L^2(\widehat{\Omega}_{\varepsilon,3})]^n}
\le & || \nabla \Phi||_{[L^2(\widehat{\Omega}_{\varepsilon,3})]^n}+C||
\nabla {\cal Q}_\varepsilon(\Phi)||_{[L^2(\widehat{\Omega}_{\varepsilon,3})]^n}||\nabla_y\chi_i||_{[L^2(Y)]^n}\cr
\le &C || \nabla\Phi||_{[L^2(\widehat{\Omega}_{\varepsilon,4})]^n} \le C\sqrt \varepsilon||\Phi||_{H^2(\Omega)}\le
C\sqrt\varepsilon||f||_{L^2(\Omega)}\cr}$$\fin
 We denote by  $\rho(x)=dist(x,\partial \Omega)$ the distance  between $x\in \Omega$ and the boundary
of $\Omega$. 
\medskip
\noindent{\bf Theorem 3.2 : }{\it  The solution $\phi^\varepsilon$ of problem $(3.1)$ verifies the following  estimates :
$$\leqalignno{ &
||\phi^\varepsilon-\Phi||_{L^2(\Omega)}\le C\varepsilon||f||_{L^2(\Omega)}, &(3.6)\cr  
& \Bigl\|\rho\Bigl(\nabla \phi^\varepsilon-\nabla  \Phi- \sum_{i=1}^n  {\cal
Q}_\varepsilon\Bigl({\partial\Phi\over\partial x_i}\Bigr)\nabla_y\chi_i\bigl({.\over\varepsilon}\bigr)\Bigr)
\Bigr\|_{[L^2(\Omega)]^n}\le C\varepsilon||f||_{L^2(\Omega)}. &(3.7)\cr}$$ The constants depend on $n$, $A$ and
$\partial\Omega$.}
\medskip
\noindent{\bf Proof : }  We put $\displaystyle\rho_\varepsilon(.)=\inf\Bigl\{{\rho(.)\over\varepsilon},1\Bigr\}$.
\medskip
\noindent{\bf Step one.}  Let $U\in H^1_{\Gamma_0}(\Omega)\cap H^2(\Omega)$. In problem $(3.1)$ we take  the test
function $U$, then by unfolding we transform the equality we have obtained. Thanks to $(2.2)$, $(3.4)$  and thanks to the
corollary of Proposition 3.1  of  [5], we have
$$\eqalign{
\Bigl|\int_\Omega A\bigl(\bigl\{{.\over\varepsilon}\bigr\}\bigr)\nabla\phi^{\varepsilon}.\nabla U-\int_{\Omega\times Y} A{\cal
T}_\varepsilon (\nabla\phi^\varepsilon)\nabla U\Bigr|&\le \Bigl|\int_\Omega
A\bigl(\bigl\{{.\over\varepsilon}\bigr\}\bigr)\nabla\phi^{\varepsilon}.\nabla U-\int_{\Omega\times Y} {\cal
T}_\varepsilon\Bigl(A\bigl(\bigl\{{.\over\varepsilon}\bigr\}\bigr)\nabla\phi^{\varepsilon}.\nabla U\Bigr)\Bigr|\cr
&+\Bigl|\int_{\Omega\times Y}A{\cal T}_\varepsilon(\nabla\phi^\varepsilon)\bigl\{{\cal
T}_\varepsilon(\nabla U)-\nabla U\bigr\}\Bigr|\cr   &\le
C\bigl\{\sqrt\varepsilon||\nabla\phi^\varepsilon||_{[L^2(\widehat{\Omega}_{\varepsilon,1})]^n}+
\varepsilon||\nabla\phi^\varepsilon||_{[L^2(\Omega)]^n}\bigr\}||\nabla U||_{[H^1(\Omega)]^n}\cr &\le
C\varepsilon||f||_{L^2(\Omega)}||\nabla U||_{[H^1(\Omega)]^n}\cr}$$ We apply now Theorem 2.3  to the function
$\phi^\varepsilon$. There exists $\widehat{\phi}^\varepsilon\in H^1_{per}( Y ; L^2(\Omega))$ such that
$$||{\cal T}_\varepsilon(\nabla_x \phi^\varepsilon)-\nabla_x \phi-\nabla_y\widehat{\phi}_\varepsilon||_{ [L^2(Y ;
(H^1(\Omega))^{'})]^n}
\le C\varepsilon||f||_{L^2(\Omega)}\leqno(3.8)$$ since from Lemma 3.1  we have $||\nabla
\phi^\varepsilon||_{[L^2(\widehat{\Omega}_{\varepsilon,3})]^n}\le C\sqrt\varepsilon||f||_{L^2(\Omega)}$. From the above
estimates and from  $(3.1)$ we obtain
$$\Bigl|\int_\Omega f\, U-\int_{\Omega\times Y}A\bigl(\nabla_x \phi^\varepsilon+\nabla_y\widehat{\phi}_\varepsilon)
\nabla_x U\Bigr| \le C\varepsilon||f||_{L^2(\Omega)}||\nabla U||_{[H^1(\Omega)]^n}\leqno(3.9)$$ Now let
$\overline{\chi}_i\in H^1_{per}(Y)$, $i\in\{1,\ldots, n\}$, be the solution of the variationnal problem
$$\int_Y A  \nabla_y \theta \nabla_y\bigl(\overline{\chi}_i+y_i\bigr)=0\qquad \forall  \theta\in H^1_{per}(Y)\leqno(3.10)$$
If matrix  $A$ is symetric $\overline{\chi}_i=\chi_i$,   $\chi_i$ are the correctors.

\noindent In problem $(3.1)$ let us take  the test function $u_\varepsilon(x)=\displaystyle
\varepsilon\rho_\varepsilon(x)\sum_{i=1}^n {\cal Q}_{\varepsilon}({\partial U\over \partial x_i})(x)\overline{\chi}_i\bigl(
{x\over\varepsilon} \bigr)$. We have multiplied by  $\rho_\varepsilon$ so that the test function $u_\varepsilon$  belongs
to $H^1_0(\Omega)$. We immediately verify the inequalities ($i\in\{1,\ldots, n\}$)
$$\eqalign{
\Big|\int_\Omega  A\bigl(\bigl\{{.\over\varepsilon}\bigr\}\bigr)\nabla \phi^{\varepsilon}\nabla u_\varepsilon\Big|=\Bigl|\int_\Omega f u_\varepsilon\Bigr|&\le
C\varepsilon||f||_{L^2(\Omega)}||\nabla U||_{[L^2(\Omega)]^n}\cr
\Bigl|\int_\Omega \varepsilon A\bigl(\bigl\{{.\over\varepsilon}\bigr\}\bigr)\nabla \phi^{\varepsilon}\nabla  \rho_\varepsilon
\,{\cal Q}_{\varepsilon}({\partial U\over \partial x_i})\overline{\chi}_i\bigl( {.\over\varepsilon}\bigr)\Bigr|&\le C\sqrt
\varepsilon||\nabla \phi^\varepsilon||_{ [L^2(\widehat{\Omega}_{\varepsilon,1})]^n}||\nabla U||_{[H^1(\Omega)]^n}\cr
\Bigl|\int_\Omega  \varepsilon\rho_\varepsilon A\bigl(\bigl\{{.\over\varepsilon}\bigr\}\bigr)\nabla \phi^{\varepsilon} \,\nabla
{\cal Q}_{\varepsilon}({\partial U\over \partial x_i})\overline{\chi}_i\bigl( {.\over\varepsilon}\bigr)\Bigr|&\le C
\varepsilon||\nabla \phi^\varepsilon||_{[L^2(\Omega )]^n}||\nabla U||_{[H^1(\Omega)]^n}\cr
\Bigl|\int_\Omega  (1-\rho_\varepsilon) A\bigl(\bigl\{{.\over\varepsilon}\bigr\}\bigr)\nabla
\phi^{\varepsilon} \,  {\cal Q}_{\varepsilon}({\partial U\over\partial x_i})\nabla_y\overline{\chi}_i\bigl( {.\over\varepsilon}
\bigr)\Bigr|&\le C\sqrt\varepsilon||\nabla \phi^\varepsilon||_{[L^2(\widehat{\Omega}_{\varepsilon,1})]^n}
||\nabla U||_{[H^1(\Omega)]^n}\cr}$$ From these estimates, from $(3.5)$  and the  corollary of Proposition  3.1 in [5]  we obtain 
$$\eqalign{ &\Bigl|\int_\Omega  A\bigl(\bigl\{{.\over\varepsilon}\bigr\}\bigr)\nabla \phi^{\varepsilon} \,\sum_{i=1}^n{\cal
Q}_{\varepsilon}({\partial
 U\over \partial x_i})\nabla_y\overline{\chi}_i\bigl( {.\over\varepsilon} \bigr)\Bigr| \le C
\varepsilon||f||_{ L^2( \Omega)}||\nabla U||_{[H^1(\Omega)]^n}\cr
\Longrightarrow \quad &\Bigl|\int_\Omega  A\bigl(\bigl\{{.\over\varepsilon}\bigr\}\bigr)\nabla \phi^{\varepsilon} \,\sum_{i=1}^n
M^\varepsilon_Y( {\partial U\over \partial x_i})\nabla_y\overline{\chi}_i\bigl( {.\over\varepsilon} \bigr)\Bigr| \le C
\varepsilon||f||_{ L^2( \Omega)}||\nabla U||_{[H^1(\Omega)]^n}\cr}$$  By unfolding we transform the left handside integral of
the above second inequality. From $(2.2)$ and $(3.5)$ we have 
$$\eqalign{ &\Bigl|\int_\Omega  A\bigl(\bigl\{{.\over\varepsilon}\bigr\}\bigr)\nabla
\phi^{\varepsilon}\,\sum_{i=1}^nM^\varepsilon_Y({\partial
 U\over \partial x_i})\nabla_y\overline{\chi}_i\bigl( {.\over\varepsilon} \bigr)-\int_{\Omega\times Y}  {\cal
T}_\varepsilon\Bigl( A\bigl(\bigl\{{.\over\varepsilon}\bigr\}\bigr)\nabla \phi^{\varepsilon}\,\sum_{i=1}^nM^\varepsilon_Y(
{\partial U\over
\partial x_i})\nabla_y\overline{\chi}_i\bigl( {.\over\varepsilon} \bigr)\Bigr)\Bigr|\cr
\le &C\sqrt\varepsilon||\nabla \phi^\varepsilon||_{[L^2(\widehat{\Omega}_{\varepsilon,1})]^n}||\nabla U||_{[H^1(
\Omega)]^n}\le C\varepsilon||f||_{L^2( \Omega)}||\nabla U||_{[H^1(\Omega)]^n}\cr}$$ We reintroduce the partial
derivatives of $U$. As a result we have 
$$\Bigl|\int_{\Omega\times Y}  A{\cal T}_\varepsilon(\nabla_x \phi^{\varepsilon}) \, \sum_{i=1}^n{\partial U\over \partial
x_i}\nabla_y \overline{\chi}_i\Bigr|\le C\varepsilon||f||_{L^2(\Omega)}||\nabla U||_{[H^1(\Omega)]^n}$$ We replace
${\cal T}_\varepsilon(\nabla_x \phi^{\varepsilon})$ by $\nabla_x \phi+\nabla_y\widehat{\phi}_\varepsilon$ thanks to
$(3.8)$, which gives us
$$\Bigl|\int_{\Omega\times Y}  A\bigl(\nabla_x \phi^{\varepsilon}+\nabla_y\widehat{\phi}^\varepsilon\bigr)
\,\nabla_y\Bigl(\sum_{i=1}^n{\partial U\over \partial x_i}\overline{\chi}_i\Bigr)\Bigr|\le C\varepsilon ||f||_{L^2(\Omega)}
||\nabla U||_{[H^1(\Omega)]^n}$$ From the definition of the correctors $\chi_i$ we obtain
$\displaystyle
\int_{\Omega
\times Y}  A\bigl(\nabla_x \phi^\varepsilon+\sum_{i=1}^n{\partial\phi^\varepsilon\over\partial x_i}\nabla_y
\chi_i\bigr)\,\nabla_y\Bigl(\sum_{j=1}^n{\partial U\over \partial x_j}\overline{\chi}_j\Bigr)=0$, we substract it from the left
handside of the above inequality and thanks to  $(3.10)$ we deduce
$$\Bigl|\int_{\Omega\times Y}  A\nabla_y\bigl(\widehat{\phi}^\varepsilon-\sum_{i=1}^n{\partial\phi^\varepsilon\over\partial
x_i}\chi_i\bigr)\,\nabla_x U\Bigr|\le C\varepsilon ||f||_{L^2(\Omega)}||\nabla U||_{[H^1(\Omega)]^n}$$  and then from
 $(3.9)$ we obtain
$$\Bigl|\int_{\Omega}{\cal A}\bigl(\nabla \phi^\varepsilon-\nabla \Phi\bigr)\nabla U\Bigr| \le C\varepsilon
||f||_{L^2(\Omega)}||\nabla U||_{[H^1(\Omega)]^n}\qquad \forall  U\in H^1_{\Gamma_0}(\Omega)\cap
H^2(\Omega)\leqno(3.11)$$ where ${\cal A}$ is the matrix of the homogenized problem.
\medskip
\noindent  Let $U_\varepsilon\in H^1_{\Gamma_0}(\Omega)$ be the solution of the variationnal problem
$$\int_\Omega {\cal A}\nabla v\nabla U_\varepsilon=\int_\Omega(\phi^\varepsilon-\Phi)v\qquad\forall
v\in H^1_{\Gamma_0}(\Omega)\leqno(3.12)$$ The boundary of $\Omega$ is of class ${\cal C}^{1,1}$ and we have
the homogeneous  Dirichlet  or homogeneous  Neumann  limits conditions. As a result we have $U_\varepsilon$
belonging to $H^1_{\Gamma_0}(\Omega)\cap H^2(\Omega)$. Moreover it verifies the estimate
$$||U_\varepsilon||_{H^2(\Omega)}\le C||\phi^\varepsilon-\Phi||_{L^2(\Omega)}$$  In $(3.12)$ we take
$v=\phi^\varepsilon-\Phi$ to obtain the estimate  of the $L^2$ norm of  $\phi^\varepsilon-\Phi$ thanks to
$(3.11)$.
\medskip
\noindent{\bf Step two.}  Now we prove the estimate  $(3.7)$ of the theorem. 
\medskip
\noindent Let $U$ be in $ H^1_{\Gamma_0}(\Omega)$. From Theorem 3.4 in  [5] there exists
$\widehat{u}^\varepsilon\in H^1_{per}(Y;L^2(\Omega))$ such that
$$||{\cal T}_\varepsilon(\nabla U)-\nabla U-\nabla_y\widehat{u}_\varepsilon||_{ [L^2(Y ;  H^{-1}(\Omega))]^n}
\le C\varepsilon||\nabla  U||_{[L^2(\Omega)]^n}\leqno(3.13)$$ In problem  $(3.1)$ we take the test function   $\rho U$ and
in problem  $(3.2)$ the couple of test functions $(\rho U, \rho\widehat{u}^\varepsilon)$. We obtain
$$\left\{\eqalign{
\int_{\Omega}f\,\rho  U&=\int_{\Omega }A\bigl(\bigl\{{.\over \varepsilon}\bigr\}\bigr)\,\rho \nabla  \phi^\varepsilon.\nabla U +
\int_{\Omega }U A\bigl(\bigl\{{.\over \varepsilon}\bigr\}\bigr) \nabla \phi^\varepsilon\nabla\rho\cr
\int_{\Omega}f\,\rho U&=\int_{\Omega\times Y }A\rho \Bigl(\nabla_x  \Phi+\sum_{i=1}^n{\partial\Phi\over\partial
x_i}\nabla_y\chi_i\Bigr)\bigl(\nabla_x U+\nabla_y\widehat{u}^\varepsilon\bigr) \cr &+ \int_{\Omega\times Y }U A
\Bigl(\nabla_x  \Phi+\sum_{i=1}^n{\partial\Phi\over\partial x_i}\nabla_y\chi_i\Bigr)\nabla_x\rho\cr}\right.\leqno(3.14)$$ The
 solution $\Phi$ of homogenized problem $(3.3.i)$ belongs to
$H^2(\Omega)$ and verifies $||\Phi||_{H^2(\Omega)}\le C||f||_{L^2(\Omega)}$. Hence the function
$\rho\nabla \Phi$ belongs to  $[H^1_0(\Omega)]^n$. From $(3.13)$ we have
$$\Bigl|\int_{\Omega\times Y}A\rho \bigl(\nabla_x \Phi+\sum_{i=1}^n{\partial\Phi\over\partial x_i}\nabla_y\chi_i\bigr)
\bigl({\cal T}_\varepsilon(\nabla_x U)-\nabla_x U-\nabla_y\widehat{u}_\varepsilon\bigr)\Bigr|\le
C\varepsilon||f||_{L^2(\Omega)}||\nabla U||_{[L^2(\Omega)]^n}$$ Now we introduce the discrete functions
$M^\varepsilon_Y(\nabla \Phi)$, $M^\varepsilon_Y({\partial\Phi\over \partial x_i})$, $M^\varepsilon_Y(U)$,
$M^\varepsilon_Y(\rho)$, $M^\varepsilon_Y(\nabla\rho)$  to replace $\nabla\Phi$, $ {\partial\Phi\over
\partial x_i}$, $U$, $\rho$, $\nabla\rho$  thanks to the estimate  of Proposition 3.1 of [5]). We use  $(2.2)$ to transform  the
integrals over $\Omega\times Y$ in integrals over $\Omega$ by inverse unfolding. Then we replace the discrete functions by 
$\nabla\Phi$, ${\cal Q}_\varepsilon( {\partial\Phi\over \partial x_i})$, $U$, $\rho$, $\nabla\rho$ and to conclude we add 
the partial derivatives   missing  in the gradient of $\Phi+\varepsilon\sum_{i=1}^n{\cal Q}_\varepsilon\Bigl(
{\partial\Phi\over\partial x_i}\Bigr) \chi_i\bigl({.\over \varepsilon} \bigr)$   (for more details see the proof of Proposition 4.3 in
[5]). We obtain
$$\eqalign{
\Bigl|\int_\Omega f\,\rho U-&\int_\Omega  A\bigl(\bigl\{{.\over \varepsilon}\bigr\}\bigr)\rho \nabla \Bigl(
\Phi+\varepsilon\sum_{i=1}^n{\cal Q}_\varepsilon\Bigl( {\partial\Phi\over\partial x_i}\Bigr) \chi_i\bigl({.\over
\varepsilon} \bigr)\Bigr) \nabla U\cr   - &\int_{\Omega }U A\bigl(\bigl\{{.\over \varepsilon}\bigr\}\bigr) \nabla
\Bigl(\Phi+\varepsilon\sum_{i=1}^n {\cal Q}_\varepsilon\Bigl({\partial\Phi\over\partial x_i}\Bigr) \chi_i\bigl( {.\over
\varepsilon} \bigr)\Bigr)\nabla\rho\Bigr|\le C\varepsilon||f||_{L^2(\Omega)}|| U||_{H^1(\Omega) }\cr}$$ The first equality of
 $(3.14)$ and the above inequality give us
$$\eqalign{
\Bigl|&\int_\Omega  A\bigl(\bigl\{{.\over \varepsilon}\bigr\}\bigr)\rho \nabla \Bigl(\phi^\varepsilon-
\Phi-\varepsilon\sum_{i=1}^n{\cal Q}_\varepsilon\Bigl( {\partial\Phi\over\partial x_i}\Bigr) \chi_i\bigl( {.\over
\varepsilon} \bigr)\Bigr) \nabla U\cr   +&\int_{\Omega }U A\bigl(\bigl\{{.\over\varepsilon}\bigr\}\bigr)\nabla
\Bigl(\phi^\varepsilon-\Phi-\varepsilon
\sum_{i=1}^n  {\cal Q}_\varepsilon\Bigl({\partial\Phi\over\partial x_i}\Bigr)\chi_i\bigl({.\over
\varepsilon} \bigr)\Bigr)\nabla\rho\Bigr|\le C\varepsilon||f||_{L^2(\Omega)}|| U||_{H^1(\Omega) }\cr}$$ Now we choose
$\displaystyle U=\rho\Bigl(\phi^\varepsilon-\Phi-\varepsilon \sum_{i=1}^n  {\cal Q}_\varepsilon\Bigl({\partial\Phi\over\partial
x_i}\Bigr)\chi_i\bigl( {.\over \varepsilon} \bigr)\Bigr)$. From the coercivity of matrix $A$ there follows that
$$\eqalign{
&||\rho\nabla \Bigl(\phi^\varepsilon-\Phi-\varepsilon \sum_{i=1}^n  {\cal
Q}_\varepsilon\Bigl({\partial\Phi\over\partial x_i}\Bigr)\chi_i\bigl({.\over \varepsilon} \bigr)\Bigr)||^2_{[L^2(\Omega)]^n}\cr
\le C &||\rho\nabla \Bigl(\phi^\varepsilon-\Phi-\varepsilon \sum_{i=1}^n  {\cal Q}_\varepsilon\Bigl({\partial\Phi\over\partial
x_i}\Bigr)\chi_i\bigl({.\over \varepsilon} \bigr)\Bigr)||_{[L^2(\Omega)]^n}|| \phi^\varepsilon-\Phi-\varepsilon
\sum_{i=1}^n  {\cal Q}_\varepsilon\Bigl({\partial\Phi\over\partial x_i}\Bigr)\chi_i\bigl({.\over
\varepsilon} \bigr) ||_{L^2(\Omega)}\cr + C &\varepsilon||f||_{L^2(\Omega)}\Bigl\{||\rho\nabla
\Bigl(\phi^\varepsilon-\Phi-\varepsilon \sum_{i=1}^n  {\cal Q}_\varepsilon\Bigl({\partial\Phi\over\partial
x_i}\Bigr)\chi_i\bigl({.\over \varepsilon} \bigr)\Bigr)||_{[L^2(\Omega)]^n}+||
\phi^\varepsilon-\Phi-\varepsilon\sum_{i=1}^n  {\cal Q}_\varepsilon\Bigl({\partial\Phi\over\partial x_i}\Bigr)\chi_i\bigl({.\over
\varepsilon} \bigr) ||_{L^2(\Omega)}\Bigr\}\cr}$$ Thanks to $(3.6)$ we obtain an upper bound of $||\rho\nabla
\bigl(\phi^\varepsilon-\Phi-\varepsilon \sum_{i=1}^n  {\cal Q}_\varepsilon\bigl({\partial\Phi\over\partial x_i}\bigr)\chi_i ({.\over
\varepsilon})\bigr)||_{[L^2(\Omega)]^n}$. The functions ${\cal Q}_\varepsilon\bigl({\partial\Phi\over\partial x_i}\bigr)$,
$i\in\{1,\ldots, n\}$, are bounded in  $H^1(\Omega)$, the estimate $(3.7)$ immediately follows.\fin
\noindent{\bf Corollary  : } Let $\Omega^{'}$ an open set strongly included in  $\Omega$, we have
$$||\phi^\varepsilon-\Phi-\varepsilon \sum_{i=1}^n  {\cal Q}_\varepsilon\Bigl({\partial\Phi\over\partial
x_i}\Bigr)\chi_i\bigl({.\over \varepsilon} \bigr)||_{H^1(\Omega^{'})}\le C\varepsilon||f||_{L^2(\Omega)}$$
 The constant depends on $n$, $A$, $\Omega^{'}$ and $\partial\Omega$.\fin
\noindent{\bf 3.2 Second case : Lipschitz boundary}
\medskip
\noindent In Theorem 4.5  of [5], $\Gamma_0$ is a  union of connected components of $\partial
\Omega$ and we have  shown that   there exists $\gamma$ in the interval $\displaystyle\bigl]0,{1/3}\bigr]$
depending on $A$, $n$ and $\partial \Omega$ such that  the solution of problem $(3.1)$ verifies the following error estimate :
$$||\phi^\varepsilon-\Phi||_{L^2(\Omega)}+||\nabla \phi^\varepsilon-\nabla \Phi-\sum_{i=1}^n {\cal Q}_\varepsilon\Bigl(
{\partial\Phi\over\partial x_i}\Bigr)\nabla_y\chi_i\bigl({.\over\varepsilon}\bigr)||_{[L^2(\Omega)]^n}
\le C\varepsilon^\gamma ||f||_{L^2(\Omega)}\leqno(3.15)$$ The constant depends on $n$, $A$ and $\partial\Omega$. 
\medskip
\noindent{\it  In the sequel of this paragraph we suppose that  

$\bullet$ the open set $\Omega$ is a bounded domain in $\R^2$ of polygonal ($n=2$) or polyhedral $(n=3)$ boundary,

$\bullet$ $\Omega$ is on one side only of its boundary, 

$\bullet$  $\Gamma_0$ is the union of some edges ($n=2$) or some faces $(n=3)$ of $\partial\Omega$,

$\bullet$  if $\Gamma_0\not=\partial\Omega$  the homogenized matrix ${\cal A}$ is symetric. }
\medskip
\noindent We know (see [6]) that for any $g\in L^2(\Omega)$ the solution of the variationnal problem 
$$U\in H^1_{\Gamma_0}(\Omega),\qquad \int_\Omega\nabla U\nabla\phi=\int_\Omega g\phi\qquad
\forall\phi\in H^1_{\Gamma_0}(\Omega)\leqno(3.16)$$ belongs to $H^{1+s}(\Omega)$ for an $s$ belonging to
$]1/2, 1[$ ($s=1$ if the domain is convex)  depending only on  $\partial\Omega$ and on the  chosen limits conditions and 
verifies the estimate 
$$||\nabla U||_{s,\Omega}\le C||g||_{L^2(\Omega)}$$
Under  a non singular linear transformation  the variationnal problem $(3.3)$ becomes $(3.16)$. It is posed in a domain which is
of the same kind as $\Omega$. Hence, the   solution $\Phi$ of the homogenized problem $(3.3)$ belongs to 
$H^{1+s}(\Omega)$ for an $s$ belonging to $]1/2, 1[$ ($s=1$ if the domain is convex)  depending only on  $\partial\Omega$,
on ${\cal A}$ and on the  chosen limits conditions and verifies the estimate 
$$||\nabla\Phi||_{s,\Omega}\le C||f||_{L^2(\Omega)}$$
\noindent{\bf Theorem 3.3 : }{\it  The solution $\phi^\varepsilon$ of problem $(3.1)$ verifies
$$\left\{\eqalign{ 
&\Bigl\|\nabla  \phi^\varepsilon-\nabla \Phi-  \sum_{i=1}^n  {\cal Q}_\varepsilon\Bigl({\partial\Phi\over\partial x_i}\Bigr)
\nabla_y\chi_i\bigl({.\over\varepsilon}\bigr)\Bigr\|_{[L^2(\Omega)]^n} \le C \varepsilon^{s/2} ||f||_{L^2(\Omega)},
\cr  &||\phi^\varepsilon-\Phi||_{L^2(\Omega)}+\Bigl\|\rho \Bigl(\nabla\phi^\varepsilon-\nabla\Phi- \sum_{i=1}^n  {\cal
Q}_\varepsilon\Bigl({\partial\Phi\over\partial x_i}\Bigr)\nabla_y\chi_i\bigl({.\over\varepsilon}\bigr)\Bigr)
\Bigr\|_{[L^2(\Omega)]^n}  \le C\varepsilon^s||f||_{L^2(\Omega)}. \cr}\right.\leqno(3.17)$$ The constants depend on $n$,
$A$ and $\partial\Omega$.}
\medskip
\noindent{\bf Proof : }

\noindent{\bf Step one.} As in Proposition 4.3 of [5], we show  that if $(\Phi,\widehat{\phi})$ is the solution of problem
$(3.2)$, then 
$\displaystyle \Phi+\sum_{i=1}^n\varepsilon\rho_{\varepsilon}{\cal Q}_\varepsilon\Bigl({\partial\Phi\over \partial
x_i}\Bigr)\chi_i\bigl({.\over\varepsilon}\bigr)$ is an approximate solution of problem $(3.1)$.  The function $\Phi$ is the
solution of the homogenized problem $(3.3)$. 

Let $\Psi\in H^1_{\Gamma_0}(\Omega)$. Thanks to Theorem 2.6, there exists $\widehat{\psi}^\varepsilon\in H^1_{per}(Y;
L^2(\Omega))$ verifying the estimates $(2.15)$. We take $(\Psi,\widehat{\psi}^\varepsilon)$ as test-function in the unfolded
problem $(3.2)$. Since $\nabla \Phi$ belongs to $ [H^s(\Omega)]^n$ and
$||\nabla\Phi||_{s,\Omega}\le C ||f||_{L^2(\Omega)}$, we obtain
$$\Bigl|\int_\Omega f \Psi-{1\over |Y|}\int_{\Omega\times Y}A(y)\Bigl(\nabla_x\Phi(x)+\sum_{i=1}^n{\partial\Phi\over \partial
x_i}(x)\nabla_y\chi_i(y)\Bigr){\cal T}_\varepsilon\bigl(\nabla_x\Psi\bigr)\Bigr|\le
C\varepsilon^{s/2}||f||_{L^2(\Omega)}||\Psi||_{H^1(\Omega)}$$
We replace  $\displaystyle{\partial \Phi\over \partial x_i}$ by  $\displaystyle{\cal Q}_\varepsilon\bigl({\partial \Phi\over \partial
x_i}\bigr)$  and then, the following part of the proof is exactly the same as the proof of
Proposition 4.3 in [5] because, thanks to Lemma 2.4 we have
$$\left\{\eqalign{ 
&||\nabla\Phi-{\cal Q}_\varepsilon(\nabla\Phi)||_{[L^2( \Omega )]^n}\le C\varepsilon^s  ||f||_{L^2(\Omega)},\cr 
&||{\cal Q}_\varepsilon(\nabla\Phi)||_{[L^2(\widehat{\Omega}_{\varepsilon,3})]^n}\le
C\varepsilon^{s/2}  ||f||_{L^2(\Omega)},\cr  &||{\cal Q}_\varepsilon(\nabla\Phi)||_{[L^2( \Omega)]^n}\le C
||f||_{L^2(\Omega)}\qquad ||{\cal Q}_\varepsilon(\nabla\Phi)||_{[H^1( \Omega)]^n}\le C\varepsilon^{s-1} 
||f||_{L^2(\Omega)}.\cr}\right.\leqno(3.18)$$  Hence we obtain the first inequality of  $(3.17)$. 

\noindent{\bf Step two. } We now use the first inequality of $(3.17)$ and again the estimates of Lemma 2.4 and  as
 in Lemma 3.1 we prove the following upper bound of the $L^2$ norm of gradient 
$\phi^\varepsilon$ in the neighbourhood of $\Omega$ :
$$  ||\nabla \phi^\varepsilon||_{[L^2(\widehat{\Omega}_{\varepsilon,3})]^n}\le C \varepsilon^{s/2}
||f||_{L^2(\Omega)}\leqno(3.19)$$  The constant depends on $n$,  $A$ and
$\partial\Omega$.

\noindent{\bf Step three. } Let $U$ be in  $H^1_{\Gamma_0}(\Omega)\cap H^{1+s}(\Omega)$. In problem $(3.1)$ we take 
the test function $U$, then by unfolding we transform the equality we have obtained. Thanks to $(2.2)$, $(3.19)$  and thanks
to the corollary of Lemma 2.4, we have
$$\eqalign{
\Bigl|\int_\Omega A\bigl(\bigl\{{.\over\varepsilon}\bigr\}\bigr)\nabla\phi^{\varepsilon}.\nabla U-\int_{\Omega\times Y} 
A{\cal T}_\varepsilon (\nabla\phi^\varepsilon)\nabla U\Bigr|&\le \Bigl|\int_\Omega
A\bigl(\bigl\{{.\over\varepsilon}\bigr\}\bigr)\nabla\phi^{\varepsilon}.\nabla U-\int_{\Omega\times Y} {\cal
T}_\varepsilon\Bigl(A\bigl(\bigl\{{.\over\varepsilon}\bigr\}\bigr)\nabla\phi^{\varepsilon}.\nabla U\Bigr)\Bigr|\cr
&+\Bigl|\int_{\Omega\times Y}A{\cal T}_\varepsilon(\nabla\phi^\varepsilon)\bigl\{{\cal
T}_\varepsilon(\nabla U)-\nabla U\bigr\}\Bigr|\cr    
&\le C\bigl\{\varepsilon^{s/2}||\nabla\phi^\varepsilon||_{[L^2(\widehat{\Omega}_{\varepsilon,1})]^n}+
\varepsilon^s||\nabla\phi^\varepsilon||_{[L^2(\Omega)]^n}\bigr\}||\nabla U||_{s,\Omega}\cr  &\le
C\varepsilon^s||f||_{L^2(\Omega)}||\nabla U|||_{s,\Omega}\cr}$$ We now apply Theorem 2.6  to the function
$\phi^\varepsilon$. There exists $\widehat{\phi}^\varepsilon\in H^1_{per}( Y ; L^2(\Omega))$ such that
$$||{\cal T}_\varepsilon(\nabla_x \phi^\varepsilon)-\nabla_x \phi-\nabla_y\widehat{\phi}_\varepsilon||_{ [L^2(Y ;
(H^s(\Omega))^{'})]^n} \le C\varepsilon^s||f||_{L^2(\Omega)}\leqno(3.20)$$ We go on as in step 1 of Theorem 3.2 to obtain
$$\Bigl|\int_{\Omega}{\cal A}\bigl(\nabla \phi^\varepsilon-\nabla \Phi\bigr)\nabla U\Bigr| \le C\varepsilon^s
||f||_{L^2(\Omega)}||\nabla U||_{s,\Omega}\qquad \forall U\in H^1_{\Gamma_0}(\Omega)\cap
H^{1+s}(\Omega)\leqno(3.21)$$
\noindent  Let $U_\varepsilon$ be  the solution of the variationnal problem
$$U_\varepsilon\in  H^1_{\Gamma_0}(\Omega),\qquad\qquad \int_\Omega {\cal A}\nabla v\nabla
U_\varepsilon=\int_\Omega(\phi^\varepsilon-\Phi)v\qquad\forall v\in H^1_{\Gamma_0}(\Omega).\leqno(3.22)$$ The function
$U_\varepsilon$ belongs to
$H^1_{\Gamma_0}(\Omega)\cap H^{1+s}(\Omega)$. Moreover we have
$$||\nabla U_\varepsilon||_{s,\Omega}\le C||\phi^\varepsilon-\Phi||_{L^2(\Omega)}$$ We take $v=\phi^\varepsilon-\Phi$ in 
$(3.22)$ and  thanks to $(3.21)$   we obtain the estimate of the $L^2$ norm of $\phi^\varepsilon-\Phi$.
\medskip
\noindent{\bf Step four.}  We now prove the upper bound of $\displaystyle \rho\Bigl(\nabla
\phi^\varepsilon-\nabla\Phi- \sum_{i=1}^n  {\cal Q}_\varepsilon\Bigl({\partial\Phi\over\partial x_i}\Bigr)
\nabla_y\chi_i\bigl({.\over \varepsilon} \bigr)\Bigr)$. 
\medskip
\noindent We take a test function in  $U\in H^1_{\Gamma_0}(\Omega)$ and as in step 2 of Theorem 2.5 we decompose the 
unfolded of its  gradient thanks to Theorem 3.4 of  [5]. In $(3.1)$  we take
$\rho U$ as test function and in $(3.2)$ we take   $(\rho U, \rho\widehat{u}^\varepsilon)$ as couple of test functions.
We obtain both equalities  $(3.14)$. In the first line of the second equality of $(3.14)$  we replace $\nabla \Phi$ and
${\partial\Phi\over \partial x_i}$ by ${\cal Q}_\varepsilon(\nabla\Phi)$ and ${\cal Q}_\varepsilon({\partial\Phi\over \partial
x_i})$.  Thanks to $(3.18)$ we have
$$\left\{\eqalign{
\Bigl| &\int_{\Omega}f\,\rho U -\int_{\Omega\times Y }A\rho \Bigl({\cal Q}_\varepsilon(\nabla_x  \Phi)+
\sum_{i=1}^n{\cal Q}_\varepsilon({\partial\Phi\over\partial x_i})\nabla_y\chi_i\Bigr)\bigl(\nabla_x
U+\nabla_y\widehat{u}^\varepsilon\bigr) \cr &+ \int_{\Omega\times Y }U A \Bigl(\nabla_x 
\Phi+\sum_{i=1}^n{\partial\Phi\over\partial x_i}\nabla_y\chi_i\Bigr)\nabla_x\rho\Bigr|\le
C\varepsilon^s||f||_{L^2(\Omega)}||\nabla U||_{[L^2(\Omega)]^n}\cr}\right.$$ From the belonging of
$\rho{\cal Q}_\varepsilon(\nabla\Phi)$ to $[H^1_0(\Omega)]^n$, and from $(3.18)$ and from   $(3.13)$ we deduce
$$\Bigl|\int_{\Omega\times Y}A\rho \Bigl({\cal Q}_\varepsilon(\nabla_x \Phi)+\sum_{i=1}^n{\cal
Q}_\varepsilon({\partial\Phi\over\partial x_i})\nabla_y\chi_i\Bigr)
\Bigl({\cal T}_\varepsilon(\nabla_x U)-\nabla_x U-\nabla_y\widehat{u}_\varepsilon\Bigr)\Bigr|\le
C\varepsilon^s||f||_{L^2(\Omega)}||\nabla U||_{[L^2(\Omega)]^n}$$ We go on as in step 2 of Theorem 3.2.
To conclude we use the upper bound of the $L^2$ norm of the function $\phi^\varepsilon-\Phi$ we obtained above.\fin
\noindent{\bf Corollary : } Let $\Omega^{'}$ be an open set strongly included in $\Omega$, we have
$$||\phi^\varepsilon-\Phi-\varepsilon \sum_{i=1}^n  {\cal Q}_\varepsilon\Bigl({\partial\Phi\over\partial
x_i}\Bigr)\chi_i\Bigl({.\over \varepsilon} \Bigr)||_{H^1(\Omega^{'})}\le C\varepsilon^s||f||_{L^2(\Omega)}$$
 The constant depends on $n$, $A$, $\Omega^{'}$ and $\partial\Omega$.\fin
 \noindent{\bf Remark : } If $\Omega$ is a convex domain we obtain the same estimates as in Theorem 3.2.\fin
\bigskip
 \centerline{\bf References}
\medskip
\noindent [1] A. Bensoussan, J.-L.Lions and G.Papanicolaou, Asymptotic Analysis for Periodic Structures, North Holland,
Amsterdam, 1978.

\noindent [2] D. Cioranescu, A. Damlamian et  G. Griso,  Periodic Unfolding and Homogenization. C. R. Acad. Sci. Paris, Ser. I 335
(2002), 99--104.

\noindent [3] D. Cioranescu and P. Donato,  An  Introduction to Homogenization. Oxford Lecture Series in Mathematics ans its
Applications 17, Oxford University Press, 1999.

\noindent [4]  G. Griso, Estimation d'erreur et \'eclatement en homog\'en\'eisation p\'eriodique.  C. R. Acad. Sci. Paris, Ser. I 335
(2002), 333--336.

\noindent [5]  G. Griso, Error estimate and unfolding for periodic homogenization, Asymptotic Analysis, Vol. 40, 3-4 (2004), 269-286.

\noindent [6]   Grivard, Singularities in Boundary Value Problems, Masson and Springer-Verlag, 1992.  

\noindent [7] O.A. Oleinik, A. S. Shamaev and G. A. Yosifian, Mathematical Problems in Elasticity and Homogenization, North-Holland,
Amsterdam, 1992.

\bye